\documentclass[12pt,draftclsnofoot,journal,letterpaper,onecolumn]{IEEEtran}
\usepackage{epsfig}
\usepackage{graphicx, wrapfig}
\usepackage{booktabs}
\usepackage{multirow}

\usepackage[colorlinks=true,linkcolor=blue]{hyperref}
\usepackage{color}
\usepackage{cite}
\usepackage{amsfonts,amsmath,amssymb,amsthm}
\usepackage{wrapfig}
\newtheorem{theorem}{Theorem}

\newtheorem{lemma}{Lemma}

\def\sl{\em}

\def\fnote#1{\footnote}

\def\beq{\begin{equation}}
\def\eeq{\end{equation}}

\def\bR{{\mathbf{R}}}

\def\cB{{\cal B}}

\def\cX{{\cal X}}

\newcommand{\four}{\mbox{\small$\frac{1}{4}$}}

\renewcommand\arraystretch{1}
\def\*{{{\LARGE\bf $^*$}}}

\def\Diag{\mathop{\mathrm{Diag}}}

\def\Prob{\mathrm{Prob}}

\def\Diag{\mathop{\mathrm{Diag}}}

\def\cB{{\cal B}}

\def\cL{{\cal L}}

\def\cX{{\cal X}}

\def\bE{{\mathbf{E}}}

\definecolor{MyDarkBlue}{rgb}{0,0.08,0.45}
\definecolor{MyViolet}{rgb}{0.45,0.08,0.95}
\definecolor{MyBrown}{rgb}{0.45,0.08,0}

\newcommand{\be}{\begin{eqnarray}}
\newcommand{\ee}[1]{\label{#1}\end{eqnarray}}

\newcommand{\ese}{\end{eqnarray*}}
\newcommand{\bse}{\begin{eqnarray*}}
\newcommand{\rf}[1]{(\ref{#1})}

\def\Row{{\mathrm{Row}}}

\def\Diag{{\mathrm{Diag}}}

\def\bR{{\mathbf{R}}}

\def\bE{{\mathbf{E}}}

\def\Prob{\mathrm{Prob}}

\def\qed{$\Box$}

\definecolor{MyDarkBlue}{rgb}{0,0.08,0.45}
\definecolor{MyViolet}{rgb}{0.45,0.08,0.95}
\definecolor{MyBrown}{rgb}{0.45,0.08,0}
\def\qed{\hfill$\Box$\par}

\def\cL{{\cal L}}

\def\cL{{\cal L}}
\def\rm{}

\def\ML{{\mathrm{\tiny ML}}}
\def\sML{\ML}
\def\exnmb
\def\sML{{\textrm{\tiny ML}}}
\def\VI{{\mathrm{VI}}}
\def\Ap1.the1{1.2.1}

\def\etat{\eta(\omega_{t-d}^{t-1})}
\def\etatt{\eta^T(\omega_{t-d}^{t-1})}
\def\etai{\eta(\omega_{i-d}^{i-1})}
\def\etait{\eta^T(\omega_{i-d}^{i-1})}

\usepackage[dvipsnames]{xcolor}
\colorlet{Mycolor1}{green!10!orange!90!}

\newcommand{\ai}[2]{{\color{blue}\  #2}}

\title{Convex Parameter Recovery for Interacting Marked Processes}
\author{Anatoli~Juditsky, % ~\IEEEmembership{Member,~IEEE,}
	Arkadi~Nemirovski, %~\IEEEmembership{Member,~IEEE,}
	 Liyan~Xie, %\IEEEmembership{Student Member,~IEEE,}
	  Yao~Xie%\IEEEmembership{Member,~IEEE}%
\thanks{Arkadi~Nemirovski, Liyan~Xie, and Yao~Xie are with School of Industrial and Systems Engineering, Georgia Institute of Technology,
Atlanta, GA, 30332 USA e-mail: {arkadi.nemirovski@isye.gatech.edu, lxie49@gatech.edu, yao.xie@isye.gatech.edu}.
Anatoli~Juditsky is with LJK, Universit\'{e} Grenoble Alpes, Campus de Saint-Martin-d'H\`{e}res, 38401 France email: anatoli.juditsky@univ-grenoble-alpes.fr.}%
}

\begin{document}

\maketitle

\begin{abstract}

We introduce a new general modeling approach for multivariate discrete event data with categorical interacting marks, which we refer to as marked Bernoulli processes. In the proposed model, the probability of an event of a specific category to occur in a location may be influenced by past events at this and other locations. We do not restrict interactions to be positive or decaying over time as it is commonly adopted, allowing us to capture an arbitrary shape of influence from historical events, locations, and events of different categories. In our modeling, prior knowledge is incorporated by allowing general convex constraints on model parameters. We develop two parameter estimation procedures utilizing the constrained Least Squares (LS) and Maximum Likelihood (ML) estimation, which are solved using variational inequalities with monotone operators. We discuss different applications of our approach and illustrate the performance of proposed recovery routines on synthetic examples and a real-world police dataset.

\end{abstract}

%\begin{IEEEkeywords}
%Spatial-temporal Bernoulli model, Variational inequality.
%\end{IEEEkeywords}

% \IEEEpeerreviewmaketitle

%\begin{center}
%{\small \it Dedicated to people fighting with the novel coronavirus.}    \end{center}

\section{Introduction}

Discrete events are a type of sequential data, where each data point is a tuple consisting of event time, location, and possibly category. Such event data is ubiquitous in modern applications, such as police data \cite{mohler2011self}, electronic health records \cite{johnson2016mimic}, and social network data \cite{stomakhin2011reconstruction, lai2016topic}. In modeling discrete events, we are particularly interested in estimating the interactions of events, such as triggering or inhibiting effects of past events on future events. For example, in crime event modeling, the triggering effect has been empirically verified; when a crime event happens, it makes future events more likely to happen in the neighborhood. Similar empirical observations have been made for other applications such as in biological neural networks, social networks \cite{zhou2013learning,li2017detecting}, financial networks \cite{embrechts2011multivariate}, and spatio-temporal epidemiological processes \cite{kuperman2001small}.

A popular model for capturing {\it interactions} between discrete events is the so-called Hawkes processes \cite{hawkes1971spectra,hawkes1971point,hawkes1974cluster,Reinhart2017}. The Hawkes process is a type of mutually-exciting non-homogeneous point process with intensity function consisting of a deterministic part and a stochastic part depending on the past event. The stochastic part of the intensity function can capture the interactions of past events {and} the current event, and it may be parameterized in different ways. In a certain sense, Hawkes processes may be viewed as a point process analog to classical autoregression in time series analysis. Hawkes process has received much attention since it is quite general and can conveniently model interactions.
For instance, in a network Hawkes process,\footnote{When space is discretized, the spatio-temporal point process of a grid can be modeled as a network point process.} interactions between nodes are modeled using a directed weighted graph in which direction and magnitude of edges indicate direction and strength of influence of one node on another.
Along this line, there are various generalizations that allow for other types of point process modeling, where different ``link'' functions are considered, such as self-correcting process, reactive process, and specialized process (see \cite{Reinhart2017} for an overview).

Estimating the {\it interactions} of the past events and the current event is a fundamental problem for Bernoulli processes since it reveals the underlying temporal and spatial structures and predicts future events. There has been much prior work in estimating model parameters, assuming that interactions are shift-invariant and captured through {\it kernel functions}. Furthermore, various simplifying assumptions are typically made for the kernel functions, e.g., that the spatio-temporal interactions are decoupled (e.g., \cite{zhou2013learning}), implying that the interaction kernel function is a product of the interaction over time and interaction over locations and can be estimated separately. It is often assumed that the temporal kernel function decays exponentially over time with an unknown decay rate \cite{li2017detecting}, or it is completely specified \cite{hall2016tracking}; thus, the problem focus is on estimating spatial interaction between locations. It is also commonly assumed that the interactions are positive, i.e., the interaction triggers rather than inhibit future events \cite{yuan2019multivariate}. {Such simplification, however, may impede capturing} complex interaction effects between events. For instance, negative interaction or inhibition is well known to play a major role in neuronal connectivity \cite{eichler2017graphical}. The study of more complex modeling of spatial aspects, especially jointly with discrete marks, is still in infancy.

In this paper, we present a general computational framework for estimating marked spatio-temporal processes with categorical marks. Motivated by Hawkes processes, we consider a model of a discrete-time process on a finite spatio-temporal grid, which we refer to as Bernoulli processes. A brief description of the proposed modeling is as follows. At each time $t$ a site $k$ of the grid of the {\em $M$-state Bernoulli process} can be in one of $M+1$ states -- a ground state, in which ``nothing happens,'' or an event state if an event of one of $M$ given types at every (discrete) time instant $t$ takes place at the site. We assume that the probability distribution of the events at each location at time $t$ is a (linear or nonlinear) function on the process history -- past events at different sites at times from $t-d$ to $t-1$, $d$ being the memory parameter of the process. For instance, each site of a $1$-state linear (vanilla) Bernoulli process can be in one of two states -- 0 (no event) or 1 (event takes place). From the point of view of time series, this process can be seen as a vector autoregressive process. The observations at sites of the grid at time $t$ are Bernoulli random variables with the conditional expectation (what is the same as the conditional probability of an event to take place) being a linear combination of historical states of the process at times $t-d$ to $t-1$. The linear combination coefficients are unknown process parameters. This model can be seen as a natural simplification of the continuous-time Hawkes process, where spatio-temporal cells are so small that one can ignore the chances for two or more events occurring in a cell. A notable feature of our model is that prior information on the structure of interactions is represented by general convex constraints on the parameters,\footnote{Convexity is assumed for the sake of computational tractability.} 
allowing for very general types of structures of interactions.
For instance, we can relax the nonnegativity restrictions on interaction parameters and/or avoid assumptions of monotone or exponential time decay of interactions commonly used in the literature. When the situation has a ``network component'' allowing to assume that interacting sites are pairs of neighboring nodes in a known graph, we can incorporate this information, for instance, by restricting the interaction coefficients for non-neighboring pairs of sites to be zero.

The considered model is related to information diffusion processes over continuous time, for example, nonlinear Hawkes model \cite{chen2017multivariate}, self-exciting processes
over networks (see \cite{Reinhart2017} for an overview), information diffusion networks \cite{gomez2013modeling}, and multivariate stationary Hawkes processes \cite{eichler2017graphical}. Compared to these well-known models, time and space discretization leading to the spatio-temporal Bernoulli process is a considerable simplification that, nonetheless, leads to practical estimation routines that can be used in ``real world'' scenarios.

Various approaches to parametric and nonparametric estimation of spatio-temporal processes have been proposed in the literature.
A line of work \cite{fox2016spatially,yuan2019multivariate,moller2003statistical} consider non-parametric Hawkes process estimation based on the Expectation-Maximization
(EM) algorithms and the Kernel method. Least-square estimates for link functions of continuous-time multivariate stationary Hawkes process are studied in \cite{eichler2017graphical}. There is also much work \cite{mohler2013modeling,ertekin2015reactive,pitkin2018bayesian}
considering the estimation in the Bayesian framework.
In particular, \cite{python2016bayesian} considers estimation in a Bernoulli model
similar to the one we promote in this paper using the Bayesian approach and impose prior distributions on parameters.
Several authors consider the problem of sparse model estimation for point processes, see, e.g., \cite{hansen2015lasso}, etc.

An important feature of the proposed models is that they allow for simple ``computation-friendly'' statistical inferences. Our approach to processing the resulting estimation problems is based on convex optimization, which leads to computationally efficient procedures. Our primary tools here are variational inequalities (VI) with monotone operators.\footnote{Utilizing VI's with monotone operators for statistical estimation is the main novelty in our paper; to the best of our knowledge in statistics, this approach was used only once (see paper \cite{juditsky2019signal} on Generalized Linear Models).}
Specifically, we show that the parameters of spatio-temporal models can be recovered in a computationally efficient fashion by solving inequalities of this type, both in the cases of linear models (Sections \ref{sect:prb}--\ref{multistate}) and of nonlinear models satisfying certain monotonicity restrictions (Section \ref{nonllink}). In the linear case, our approach results in the usual Least Squares estimate\footnote{This, in hindsight, is resembling but not identical to what is done in \cite{eichler2017graphical}.} which, of course, could be motivated without any references to VI's.
However, these references explain how to act in the nonlinear cases, where Least Squares, if applied, ``as is'' typically lose computational tractability. Aside from the VI-based approach, we consider the standard Maximum Likelihood estimation (Section \ref{secML estimate}). In the linear case, computing the Maximum Likelihood estimate amounts to solving a convex optimization problem and thus is computationally efficient. On the other hand, in the nonlinear case, Maximum Likelihood estimation typically becomes problematic computationally, including the cases where the VI-based approach remains computation-friendly. (A notable exception is the spatio-temporal logistic model in which the Maximum Likelihood estimation reduces to solving the convex problem.) It should be added that finite-sample theoretical results on the statistical performance of the estimates we develop do not favor Maximum Likelihood as compared with the VI-based estimation.

Finally, we also demonstrate the good performance of our method on synthetic and real data. In particular, we study a real crime dataset in Atlanta, USA, to demonstrate the promise of our methods to recover interesting structures from real-data and predict the probability of crime incidents. 

The paper is organized as follows. We discuss the LS estimation of the network Bernoulli process in Section \ref{sect1}, by introducing the model of the simplest Bernoulli process with $\{0,1\}$-valued mark in Section \ref{sect:prb}. We then derive the Least Squares estimate (which under the circumstances is what VI-estimate boils down to) for this model in Section \ref{process} and building data-driven confidence sets for the estimated parameters in Section \ref{sec:bound}. We describe the general model of the $M$-state Bernoulli process and discuss the Least Squares estimation in Section \ref{multistate}. The nonlinear modeling of the process and the corresponding {VI} estimate are presented in Section \ref{nonllink}. In Section \ref{secML estimate}, we discuss the Maximum Likelihood estimate of parameters of the general Bernoulli process. The application of the proposed approach is illustrated by various simulation examples in Section \ref{sec:sim}. Finally, Section \ref{sec:rw} shows an application of our modeling to ``real-world'' data analysis of crime events in Atlanta.

\section{Estimating parameters of spatio-temporal Bernoulli process}
%\section{Recovering parameters of spatio-temporal Bernoulli process}
\label{sect1}

Here we consider spatio-temporal Bernoulli process with discrete-time over discrete locations. Specifically, we assume that the discrete-time and location grid we deal with is fine enough so that we can neglect the possibility for more than one event to occur in a cell of the grid. We will model the interactions of these events in the grid. In Sections    \ref{process}--\ref{multistate} we develop and process {\sl linear} models; nonlinear extensions are considered in Section \ref{nonllink}.

\subsection{Single-state model}\label{sect:prb}
%~\aic{}{There is a problem with Figure 1 and page split with my latex.}
%\end{figure}
%
Define a {\sl spatio-temporal Bernoulli process with memory depth $d$} as follows. We assume the memory depth is a pre-specified hyper-parameter (e.g., it can be estimated using cross-validation as explained in Section \ref{sec:rw} when we study real data). We observe on discrete time horizon $\{t: -d+1\leq t\leq N\}$ random process as follows.
%\begin{quote}
At time $t$ we observe Boolean vector $\omega_t\in\bR^K$ with entries $\omega_{tk}\in\{0,1\}$, $1\leq k\leq K$. Here $\omega_{tk}=1$ and $\omega_{tk}=0$
mean, respectively, that
at time $t$ in location $k$ an event took/did not take place. We set
\begin{equation*}
\begin{split}
\omega^t&=\{\omega_{s k},-d+1\leq s\leq t,1\leq k\leq K\}\in\bR^{(t+d)\times K},\\
\omega_\tau^t&=\{\omega_{s k},\tau\leq s\leq t,1\leq k\leq K\}\in\bR^{(t-\tau+1)\times K}.
\end{split}
\end{equation*}
In other words, $\omega^t$ denotes all observations (at all locations) until current time $t$, and $\omega_\tau^t$ contains observations on time horizon from $\tau$ to $t$.

\begin{figure}[h!]
%\begin{wrapfigure}{r}{.4\textwidth}
%\begin{figure}[h!]
\vspace{-0.1in}
  \centering\includegraphics[width=.4\textwidth]{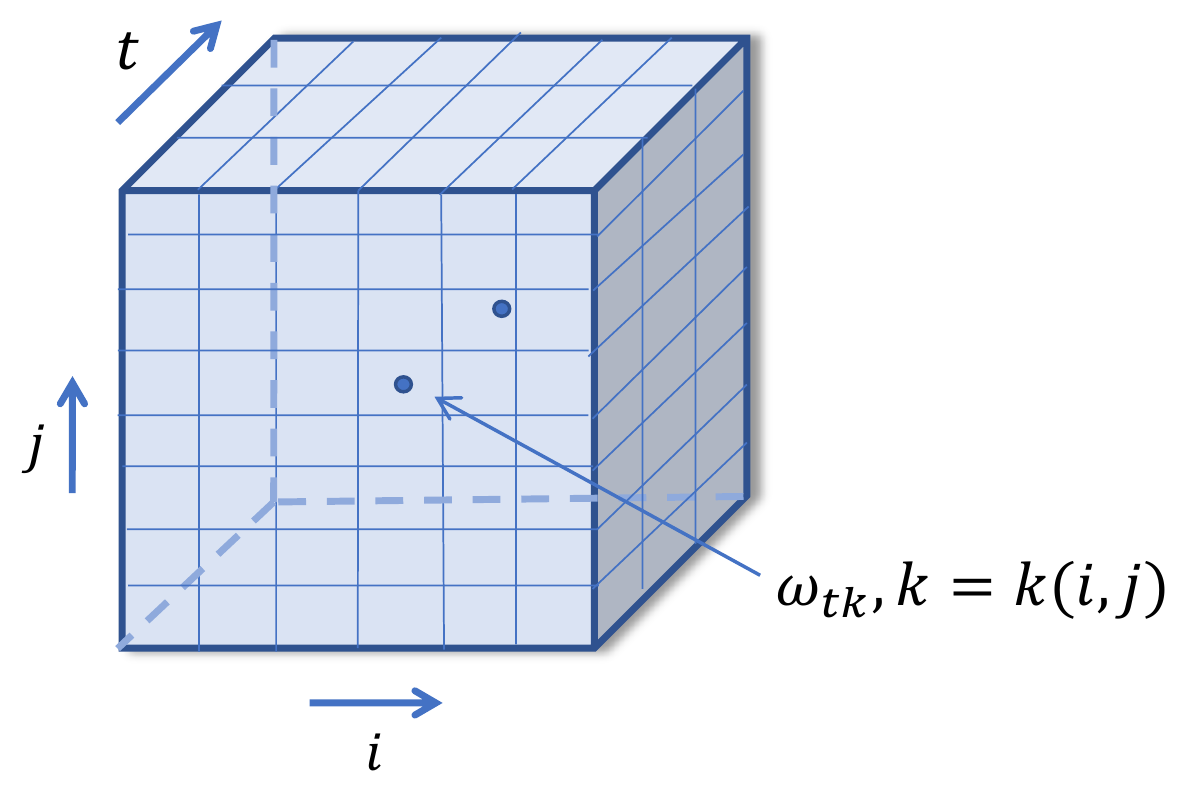}
    \vspace{-0.15in}
  \caption{Illustration of the discretized process. Observation $\omega_{tk}$, at the location of a three-dimensional spatio-temporal grid.}\label{fig:illustration}
  %\vspace{-0.15in}
%\end{wrapfigure}
\end{figure}

We  assume that for $t\geq 1$ the conditional probability of the event $\omega_{tk}=1$, given the history $\omega^{t-1}$, is specified as
\begin{equation}\label{eq1}
\beta_k+\sum_{s=1}^d \sum_{\ell=1}^K \beta^s_{k\ell}\omega_{(t-s)\ell},\,1\leq k\leq K,
\end{equation}
where $\beta=\{\beta_k,\beta^s_{k\ell}:1\leq s\leq d,1\leq k,\ell \leq K\}$ is a collection of coefficients. Here
\begin{itemize}
\item $\beta_k$ corresponds to the {\it baseline intensity} at the $k$-th location (i.e., the intrinsic probability for an event to happen at a location without the exogenous
%endogenous
influence, also called the birthrate);
\item $\beta_{k\ell}^s$ captures the magnitude of the {\it influence} of an event that occurs at time $t-s$ at the $\ell$-th location on chances for an event to happen at time $t$ in the $k$-th location; so the sum in \rf{eq1} represents the cumulative influence of past events at the $k$-th location.
\end{itemize}
%\end{quote}
\ai{
Throughout this paper, $d$ is considered as a pre-specified parameter of the procedure. In reality, it can be tuned by verifying the ``predictive abilities" of models with different $d$'s for a given dataset.
}{}
Since the probability of occurrence is between 0 and 1, we require the coefficients to satisfy
\begin{equation}\label{eq2}
\begin{array}{rcl}
0&\leq &\beta_k+\sum_{s=1}^d\sum_{\ell=1}^K\min\left[\beta^s_{k\ell},0\right],\,\,\,\forall \ k\leq K,\\
1&\geq&\beta_k+\sum_{s=1}^d\sum_{\ell=1}^K \max\left[\beta^s_{k\ell},0\right],\,\,\,\forall \ k\leq K.\\
\end{array}
\end{equation}

\begin{figure}
%\begin{wrapfigure}{r}{.4\textwidth}
\vspace{-0.1in}
\begin{center}
\includegraphics[width = 0.35\textwidth]{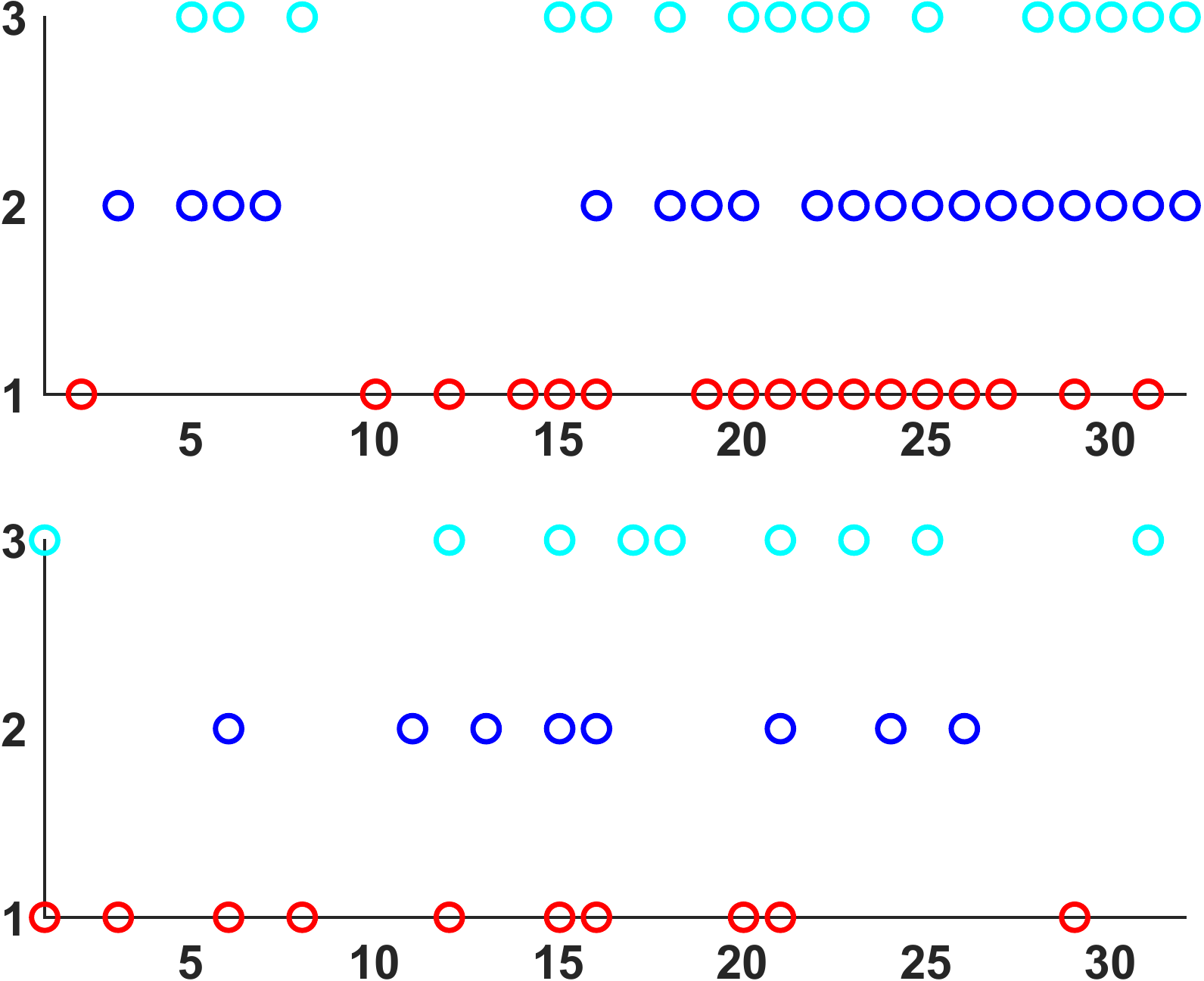}
\vspace{-0.15in}
\caption{Realizations of  spatio-temporal Bernoulli processes with memory depths 5 (top) and 0 (bottom) on time horizon $N=32$ with three locations represented with y-axis 1, 2, and 3. ``1'' events in different locations are marked by different colors.}
\label{fig:path}
\vspace{-0.15in}
\end{center}
\end{figure}

Note that constraints in \rf{eq2} allow some of the coefficients $\beta_{k\ell}^s$ to be negative, permitting the corresponding model to capture the inhibitive effect of past events.
% \begin{figure}[h!]
Fig.\,\ref{fig:path} illustrates a realization of the sample path of a simple Bernoulli process in the considered setting with different memory depths (5 for the top figure and 0 for the bottom). Note that in the bottom plot, the events are more spread out due to the memoryless nature of the process.

%$$
%\epsfxsize=200pt\epsfysize=150pt\epsffile{Bernoulli}
%$$
%\begin{quotation}
%{\tiny\noindent Realizations of spatio-temporal Bernoulli processes with memory depths 5 (top) and 0 (bottom) with 3 spatial locations (events in different locations are marked by different colors) on time horizon $N=32$. The total number of events in the realization are 30 (top) and 32 (bottom). Pay attention  to essentially more spread events on the bottom plot.}
%\end{quotation}

Our goal is to recover the collection of parameters $\beta$ using a set of observations $\omega^N$ over a time horizon $N$.

%\end{wrapfigure}

\subsection{Preliminaries on variational inequalities with monotone operators}
Variational inequalities (VI's) with monotone operators is the principal computational tool of the approach we are about to describe. We start with the related preliminaries.
A vector field $F:\,\cX\to\bR^N$ defined on a nonempty convex subset $\cX$ of $\bR^N$ is called {\sl monotone}, if $\langle F(x)-F(y), x-y \rangle \geq0$ whenever $x,y\in \cX$. When $N=1$, monotonicity means that the scalar function $F$ is nondecreasing on $\cX$; a basic example  (by far not the only useful one) of a multivariate monotone vector field is the gradient field of a differentiable convex function $f:\cX\to\bR$.
We say that $\alpha\geq0$ is a {\sl modulus of strong monotonicity} of vector field $F$, when
$$
\langle F(x)-F(y), x-y \rangle \geq\alpha\|x-y\|_2^2\,\, \forall x,y\in \cX;
$$
when $\alpha>0$, $F$ is called {\sl strongly monotone.}
A pair $(\cX,F)$ comprised of nonempty convex domain $\cX$ and monotone vector field $F$ on this domain gives rise to {\sl variational inequality}
$\VI(F,\cX)$. A {\sl weak solution} to this VI is any point $\bar{x}\in \cX$ such that
$$
\langle F(x),x-\bar{x}\rangle  \geq 0\,\,\forall x\in \cX.
$$
Whenever $F$ is strongly monotone, weak solution, if exists, is unique.\par
A {\sl strong solution} is a point $\bar{x}\in \cX$ such that
$$
\langle F(\bar{x}),x-\bar{x}\rangle  \geq0\,\,\forall x\in \cX.
$$
Every strong solution is a weak one; when $F$ is continuous on $\cX$, the inverse also is true. When $\cX$ is a convex {\sl compact} set, $\VI(F,\cX)$ always has weak solutions. When $F$ is the gradient field of a continuously differentiable convex function $f$ on $\cX$, the weak and the strong solutions to $\VI(F,\cX)$ are exactly the minimizers of $f$ on $\cX$. Finally, we should stress that variational inequalities with monotone operators are the most general ``problems with convex structure;'' under mild computability assumptions, that can be efficiently solved to a high accuracy.

\subsection{Least Squares (LS) estimation}\label{process}
As applied to the simple spatio-temporal model described in Section \ref{sect:prb} the VI-based approach we are developing boils down to the  Least Squares (LS) estimation.
Let $\kappa=K+dK^2$; we arrange all reals from the collection  $\beta$ in \eqref{eq1} into a column vector (still denoted $\beta$):
\[
\beta = [\beta_1,\ldots,\beta_K,\beta_{11}^1,\ldots,\beta_{11}^d,\beta_{1K}^1,\ldots, \beta_{1K}^d,\ldots,\beta_{KK}^1,\ldots,\beta_{KK}^d]^T \in \mathbb R^\kappa.
\]
Note that constraints (\ref{eq2}) above state that $\beta$ must reside in the polyhedral set $\cB$
given by explicit polyhedral representation.\footnote{Polyhedral representation of a set $X\subset\bR^n$ is a representation of the form
$$
X=\{x\in\bR^n: \exists w\in\bR^m: Px+Qw\leq r\},
$$
that is, representation of $X$ as a projection of the solution set of a system
of linear inequalities in the space of $(x,w)$-variables on the plane of $x$-variables.
When $X$ is polyhedrally representable, it automatically is polyhedral --- can be
represented by a finite system of linear inequalities in $x$-variables only.
This system, however, can be much larger than the one in the polyhedral representation in question,
making explicit polyhedral representations the standard descriptions of polyhedral sets in optimization.}
Assume that we are given a convex compact set
$\cX\subset\cB$ such that $\beta\in \cX$; we introduce this set to account for additional to the obvious inclusion $\beta\in\cB$
a priori information, if any, on the vector of  model's parameters.
Our model says that for $t\geq 1$,
the conditional expectation of $\omega_t$ given $\omega^{t-1}$ is $\etatt \beta$,
     \[
    \Prob_{\omega^{t-1}}\left\{\omega_{t}=1\right\}= \etatt \beta,\]
with a known to us function  $\eta(\cdot)$ which is  defined on the set of all zero-one arrays
$\omega_{t-d}^{t-1}\in \{0, 1\}^{d\times K}$ and takes values in the matrix space $\bR^{\kappa\times K}$:
\begin{equation}
\etatt  =
\begin{bmatrix}
I_K, & I_K \otimes \mbox{vec}(\omega_{t-d}^{t-1})^T
\end{bmatrix} \in \bR^{K \times \kappa},
%\begin{bmatrix}
%1 & 0 & \ldots & 0 &  \mbox{vec}(\omega_{t-d}^{t-1})^T & \ldots & 0 \\
%\vdots & \vdots & \ddots & \vdots & \vdots & \ddots & \vdots \\
%0 & 0 & \ldots & 1 & 0 & \ldots & \tilde\omega_{t-d}^{t-1} \\
%\end{bmatrix},
\label{eta_def}
\end{equation}
where $I_K$ is a $K\times K$ identity matrix, $\otimes$ denotes the standard Kronecker product,
and $\mbox{vec}(\cdot)$ vectorizes a matrix by stacking all columns.  Note that the matrix $\etat $
is Boolean and has at most one nonzero entry in every row.\footnote{Indeed, (\ref{eq1}) says that a particular entry in $\beta$, $\beta_k$ or $\beta^s_{k
\ell}$, affects at most one entry in $\etatt \beta$, namely, the $k$-th entry,
implying that each column of $\eta^T(\cdot)$ has at most one nonzero entry.}

Consider a vector field $F: \cX\to\bR^\kappa$, defined as
\[
F(x)={1\over N}\bE_{\omega^N}\left\{\sum_{t=1}^N[\etat \etatt x-\etat \omega_t]\right\}: \cX\to\bR^\kappa,
\]
where $\bE_{\omega^N}$ denotes expectation taken with respect to the distribution of $\omega^N$
(notation $\bE_{\omega^t}$ is similarly defined).
Below, all expectations and probabilities are conditional given a specific realization
of the initial fragment $\omega_{-d+1}^0$ of observations.
\par
Observe that we have
\[
\langle F(x)-F(y), x-y \rangle = {1\over N}\sum_{t=1}^N\bE_{\omega^N}\left\{(x-y)^T\etat \etatt (x-y)\right\} \geq 0, \quad \forall x,y\in\cX.
\]
Thus, the vector field $F$ is {\it monotone}.
%\footnote{A vector field $F:\,X\to\bR^N$ defined on a nonempty convex subset $X$ of $\bR^N$ is called {\sl monotone}, if $\langle F(x)-F(y), x-y \rangle \geq0$ whenever $x,y\in X$.}
Moreover, we have  $F(\beta) = 0$, since
$$
\begin{array}{rcl}
F(\beta)&=&{1\over N}\bE_{\omega^N}\left\{\sum_{t=1}^N\etat [\etatt \beta-\omega_t]\right\}\\
&=&{1\over N}
\sum_{t=1}^N\bE_{\omega^t}\left\{\etat [\etatt \beta-\omega_t]\right\}\\
&=&{1\over N}\sum_{t=1}^N\bE_{\omega^{t-1}}\left\{\etat \big[\etatt \beta-
\bE_{|\omega^{t-1}}\{\omega_t\}\big]\right\}\\
&=&{1\over N}\sum_{t=1}^N\bE_{\omega^{t-1}}\left\{\etat [\etatt \beta-\etatt \beta]\right\}=0,
\end{array}
$$
where $\bE_{|\omega^{t-1}}$ denotes the conditional expectation given $\omega^{t-1}$.
Therefore, $\beta\in\cX$ is a zero of the monotone operator $F$ and therefore it is a solution to the variational inequality $\VI[F,\cX]$.

Now consider the empirical version
\be
F_{\omega^N}(x)=\underbrace{\left[{1\over N}{\sum}_{t=1}^N\etat \etatt \right]}_{A[\omega^N]}x-
\underbrace{{1\over N}{\sum}_{t=1}^N\etat \omega_t}_{a[\omega^N]}
\ee{FF}
of vector field $F$. Note that $F_{\omega^N}(x)$ monotone and affine, and its expected value is $F(x)$ at every point $x$.

We propose to use, as an estimate of $\beta$, a weak solution to the Sample Average Approximation of $\VI[F,\cX]$, i.e., the variational inequality
$$
\hbox{find}\ z\in\cX: \langle F_{\omega^N}(w),w-z\rangle\geq 0, \quad \forall w\in\cX.\eqno{\VI[F_{\omega^N},\cX]}
$$
The monotone vector field $F_{\omega^N}(\cdot)$ is continuous (even affine), so that  weak solutions to $\VI[F_{\omega^N},\cX]$ are exactly the same as strong solutions, i.e., points $\bar{x}\in\cX$ such that $\langle F_{\omega^N}(\bar{x}),x-\bar{x}\rangle\geq0$ for all $x\in \cX$. Moreover, the empirical vector field $F_{\omega^N}(x)$ is just the gradient field of the convex quadratic function
\begin{equation}\label{Psi}
\Psi_{\omega^N}(x)={1\over 2N}\sum_{t=1}^N\|\etatt x-\omega_t\|_2^2,
\end{equation}
so that weak (same as strong) solutions to $\VI[F_{\omega^N},\cX]$ are just minimizers of this function on $\cX$. In other words, our estimate based on
solving variational inequality is an optimal solution to the Least Squares (LS) formulation:
the constrained optimization problem
\begin{equation}\label{opt1}
\min_{x\in\cX} \Psi_{\omega^N}(x)
\end{equation}
with a convex quadratic objective. Problem (\ref{opt1}), the same as a general variational inequality with a monotone operator, can be routinely and efficiently solved by convex optimization algorithms.

\subsection{Toward performance guarantees}\label{sec:bound}

Our objective in this section is to construct non-asymptotic confidence sets for parameter estimates built in the previous section.
Utilizing concentration inequalities for martingales, we can express these sets in terms of the process observations in the spirit
of results of \cite{konev1996asymptotic,liptser2000deviation,hansen2015lasso}.
\par
Observe that the vector of true parameters $\beta$  underlying our observations not only solves variational inequality $\VI[F,\cX]$, but also solves the variational inequality  $\VI[\overline{F}_{\omega^N},\cX]$,
%
%the following variational inequality
%$$
%\begin{array}{c}
%\hbox{find}\ z\in\cX: \langle \overline{F}_{\omega^N}(w),w-z\rangle\geq 0, \quad \forall w\in\cX,\\
%\end{array}
%\eqno{\VI[\overline{F}_{\omega^N},\cX]}
%$$
where
$$
\overline{F}_{\omega^N}(x)=A[\omega^N]x-\underbrace{{1\over N}{\sum}_{t=1}^N\etat \etatt \beta}_{\overline{a}[\omega^N]}
$$
with $A[\omega^N]$ defined in (\ref{FF}).

In fact, $\beta$ is just a root of $\overline{F}_{\omega^N}(x)$: $\overline{F}_{\omega^N}(\beta) = 0$.
Moreover, the monotone affine operators $F_{\omega^N}(x)$ and $\overline{F}_{\omega^N}(x)$ differ only in the value of constant term: in $F_{\omega^N}(x)$ this term is $a[\omega^N]$, and in $\overline{F}_{\omega^N}(x)$ this term is $\overline{a}[\omega^N]$. Thus, equivalently,
$\beta$ is the minimizer on $\cX$ of the quadratic form
$$
\overline{\Psi}_{\omega^N}(x):={1\over 2N}\sum_{t=1}^N\|\etatt x-\etatt \beta\|_2^2,
$$
and the functions $\Psi$ in \eqref{Psi} and $\overline{\Psi}$ above differ only in the constant terms (which do not affect the results of minimization) and in the linear terms. Moreover, the difference of the vectors of coefficients of linear terms is given by (due to $\overline{F}_{\omega^N}(\beta)=0$):
\begin{equation}\label{eq:delta}
 \Delta_F :=  {F_{\omega^N}(\beta)}-\overline{F}_{\omega^N}(\beta) = {F_{\omega^N}(\beta)} =\overline{a}[\omega^N]-a[\omega^N]=
{1\over N}\sum_{t=1}^N\underbrace{\etat [\etatt \beta-\omega_t]}_{\xi_t}.
\end{equation}
Note that this is the same as the difference of constant terms in $F_{\omega^N}(\cdot)$ and $\overline{F}_{\omega^N}(\cdot)$.

\par Concentration bounds for $F_{\omega^N}(\beta)$ %=A[\omega^N]\beta-a[\omega^N]$
can be obtained by applying general Bernstein-type inequalities for martingales.
\begin{lemma}\label{lem1}
For all $\epsilon\in(0,1)$ vector $F_{\omega^N}(\beta)=\Delta_F$ in (\ref{eq:delta}) satisfies
\begin{equation}\label{probbnd}
\Prob_{\omega^N}\left\{\|F_{\omega^N}(\beta)\|_\infty\geq \sqrt{\ln(2\kappa/\epsilon)\over 2N}+{\ln(2\kappa/\epsilon)\over 3N}\right\}\leq \epsilon.
\end{equation}
\end{lemma}
\par\noindent{\em Proof.}
Since the conditional expectation of $\omega_t$ given $\omega^{t-1}$  is $\etatt \beta$, we have $\bE_{|\omega^{t-1}}[\xi_t]=0$.
Thus, $\xi_t$ is a
{\it martingale-difference}. Also, because both $\omega_t$ and $\etatt \beta$ are vectors with nonnegative entries not exceeding 1, we have $\|\etatt \beta-\omega_t\|_\infty\leq 1$.
Besides this, $\etat $ is a Boolean matrix with at most one nonzero in every row, whence $\|\etatt z\|_\infty\leq \|z\|_\infty$ for all $z$. The bottom line is that $\|\xi_t\|_\infty\leq 1$. Furthermore, the conditional variance of components of
$\omega_t$ is bounded by $1/4$, so, applying the Azuma-Hoeffding inequality \cite{azuma1967weighted} to components $(F_{\omega^N}(\beta))_k$, $k=1,...,\kappa$, of $F_{\omega^N}(\beta)$ we conclude that
\[
\Prob_{\omega^N}\left\{|(F_{\omega^N}(\beta))_k|\geq \sqrt{x\over 2N}+{x\over 3N}\right\}\leq 2\exp\{-x\},\quad\forall 1\leq k\leq\kappa,\, x\geq0.
\]
The latter bound results in \rf{probbnd} by application of the total probability formula.\qed

A somewhat finer analysis allows to establish more precise data-driven deviation bounds for components of $ F_{\omega^N}(\beta)$.
\begin{lemma}\label{lem1002}
For all $y>1$
entries $F_{\omega^N}(\beta)_k$, $k=1,...,\kappa$, of $F_{\omega^N}(\beta)$ satisfy, with probability at least $1-2e{\big(y\big[\ln((y-1)N\big)+2\big]+2\big)}e^{-y}$,
\be
%\Prob\left\{a[\omega^N]_k-\overline\psi(a[\omega^N]_k,N;y) \leq F_{\omega^N}(\beta)_k\leq a[\omega^N]_k-\underline\psi(a[\omega^N]_k,N;y)\right\}\geq
a[\omega^N]_k-\overline\psi(a[\omega^N]_k,N;y) \leq F_{\omega^N}(\beta)_k\leq a[\omega^N]_k-\underline\psi(a[\omega^N]_k,N;y)
\ee{mart2010}
where $a[\omega^N]_k$ is the $k$-th component of $a[\omega^N]$ as in \rf{FF} and lower and upper functions $\underline\psi(\cdot)$, $\overline\psi(\cdot)$ are defined in relation \rf{psifun}, see appendix.
\end{lemma}
Proof of Lemma \ref{lem1002} is postponed till the appendix.
We are about to extract from this lemma upper bounds on the accuracy of recovered coefficients.

\subsubsection{Upper-bounding risk of recovery}
%{\ccy
Recall that our estimate $\widehat{\beta}:=\widehat{\beta}(\omega^N)$ solves the variational inequality $\VI[F_{\omega^N},\cX]$
with $F_{\omega^N}(x)=A[\omega^N]x-a[\omega^N]$, see (\ref{FF}). Note that $A[\omega^N]$  is positive semidefinite (we write $A\succeq0$, and we write $A\succ 0$ for positive definite $A$). Given $A\in\bR^{\kappa\times\kappa}$, $A\succeq0$, and $p\in[1,\infty]$, define the ``condition number''
\begin{equation} \label{theta_def}
\theta_p[A]:=\max\left\{\theta\geq0:g^TAg\geq \theta\|g\|_p^2, ~~\forall g\in\bR^\kappa\right\}.
\end{equation}
Observe that $\theta_p[A]>0$ whenever $A\succ0$, and that for $p,p'\in[1,\infty]$ one has
\begin{equation}\label{that}
g^TAg\geq {1\over 2}\left\{\theta_p[A]\|g\|_p^2+\theta_{p'}[A]\|g\|_{p'}^2 \right\} \geq \sqrt{\theta_p[A]\theta_{p'}[A]}\|g\|_p\|g\|_{p'}.
\end{equation}
The following result is immediate:
\begin{theorem}[Bounding $\ell_p$ estimation error] \label{obsini} For every $p\in[1,\infty]$ and every $\omega^N$ one has
\begin{equation}\label{onehas16}
\|\widehat{\beta}(\omega^N)-\beta\|_p\leq \|F_{\omega^N}(\beta)\|_\infty/\sqrt{\theta_p[A[\omega^N]]\theta_1[A[\omega^N]]}.
\end{equation}
As a result, for every $\epsilon\in(0,1)$, the probability of the event
\begin{equation}\label{theevent}
\|\widehat{\beta}(\omega^N)-\beta\|_p\leq \left(\theta_p[A[\omega^N]]\theta_1[A[\omega^N]]\right)^{-1}
\left(\sqrt{\ln(2\kappa/\epsilon)\over 2N}+{\ln(2\kappa/\epsilon)\over 3N}\right)
,\quad\forall p\in[1,\infty]
\end{equation}
is at least $1-\epsilon$.
\end{theorem}
\par\noindent{\em Proof.}  Let us fix $\omega^N$ and set $\widehat{\beta}=\widehat{\beta}[\omega^N]$, $A=A[\omega^N]$. Since $F_{\omega^N}(\cdot)$ is continuous and $\widehat{\beta}$ is a weak solution to $\VI[F_{\omega^N},\cX]$, $\widehat{\beta}$ is also a strong solution:
$\langle F_{\omega^N}(\widehat{\beta}), z-\widehat{\beta}\rangle\geq0$ for all $z\in \cX$; in particular,
$\langle F_{\omega^N}(\widehat{\beta}),\beta-\widehat{\beta}\rangle\geq0$. On the other hand, $F_{\omega^N}(\widehat{\beta})=F(\beta)-A(\beta-\widehat{\beta})$. As a result,
$
0\leq \langle F_{\omega^N}(\widehat{\beta}),\beta-\widehat{\beta}\rangle=\langle F_{\omega^N}(\beta)-A(\beta-\widehat{\beta}),\beta-\widehat{\beta}\rangle,
$
whence
\begin{equation}\label{eq:delta1}
(\beta-\widehat{\beta})^TA(\beta-\widehat{\beta}) \leq \langle F_{\omega^N}(\beta),\beta-\widehat{\beta}\rangle \leq \|F_{\omega^N}(\beta)\|_\infty\|\beta-\widehat{\beta}\|_1.
\end{equation}
Setting $p'=1$ in (\ref{that}), we obtain
\[
(\beta-\widehat{\beta})^TA(\beta-\widehat{\beta}) \geq \sqrt{\theta_1[A]\theta_p[A]}\|\beta-\widehat{\beta}\|_1\|\beta-\widehat{\beta}\|_p.
\] This combines with \eqref{eq:delta1} to imply (\ref{onehas16}); then
\eqref{onehas16} together with (\ref{probbnd}) imply (\ref{theevent}).
 \qed
%\par\begin{remark}
\vspace{.1in}
\par\noindent{\em Remark} [Evaluating the condition number]. To assess the upper bound \rf{theevent} one needs to compute ``condition numbers'' $\theta_p[A]$
of a positive definite matrix $A$. The computation is easy when  $p=2$, in which case $\theta_2[A]$ is the minimal eigenvalue of $A$, and when $p=\infty$:
$$\theta_\infty[A]=\min_{1\leq i \leq\kappa}\left\{x^TAx: \|x\|_\infty\leq1, x_i=1\right\}
$$
is the minimum of $\kappa$ efficiently computable quantities.  In general, $\theta_1[A]$ is difficult to compute, but this quantity admits
an efficiently computable tight within the factor $\pi/2$ lower bound. Specifically, for a symmetric positive
definite $A$, $\min_z\{z^TAz:\|z\|_1=1\}$ is the largest $r>0$ such that the ellipsoid $\{z:z^TAz\leq r\}$ is contained in the unit $\|\cdot\|_1$-ball,
or, passing to polars, the largest $r$ such that the ellipsoid $y^TA^{-1}y\leq r^{-1}$ contains the unit $\|\cdot\|_\infty$-ball.
Because of this, the definition of $\theta_1[A]$ in (\ref{theta_def}) is equivalent to $\theta_1[A]=\left[\max_{\|x\|_\infty\leq1}x^TA^{-1}x\right]^{-1}$.
It remains to note that when $Q $ is a symmetric positive semidefinite $\kappa\times\kappa$ matrix, the efficiently computable by semidefinite relaxation
upper bound on $\max_{\|x\|_\infty\leq1}x^TQx$, given by
$$
\min\limits_\lambda\left\{\sum_i\lambda_i:\lambda_i\geq0, ~\forall i; \, \Diag\{\lambda_1,...,\lambda_\kappa\}\succeq Q\right\},
$$
is tight within the factor $\pi/2$, see \cite{nesterov1998semidefinite}.
%\end{remark}
\par
Under favorable circumstances, we can expect that for large
$N$ the minimal eigenvalue of $A[\omega^N]$ will be of the order of one
with overwhelming probability implying that the lengths of the confidence intervals \rf{intervals} go to 0 as $N\to\infty$ at the rate $O(1/\sqrt{N})$.
Note, however, that inter-dependence of the ``regressors'' $\etat $ across $t$ makes it difficult to prove something along these lines.
%}%ccy
%}%ccy

\subsubsection{Estimating linear forms of $\beta$}
We can use concentration bounds of Lemmas \ref{lem1} and \ref{lem1002} to build confidence intervals
for linear functionals of  $\beta$. For instance, inequality \rf{mart2010} of Lemma \ref{lem1002} leads to the following estimation procedure of the linear form $e(\beta)=e^T\beta$, $e \in \mathbb R^\kappa$.
%, and recall that
%$$
%\begin{array}{c}
%F_{\omega^N}(w)=A[\omega^N]w - a[\omega^N],\\
%\begin{array}{rcl}
%A[\omega^N]&=&{1\over N}\sum_{t=1}^N\etat \etatt ,\\
%a[\omega^N]&=&{1\over N}\sum_{t=1}^N\etatt \omega_t\\
%\end{array}
%\\
%\end{array}
%$$
Given $y>1$, consider the pair of optimization problems
\begin{equation}\label{convopt}
\begin{array}{rcl}
\underline{e}[\omega^N,y]&=&\min\limits_x\left\{e^Tx:\begin{array}{l}\,x\in\cX,\\\underline\psi(a[\omega^N]_k,N;y)\leq (A[\omega^N]x)_k\leq\overline\psi(a[\omega^N]_k,N;y),\;k=1,...,\kappa,\end{array}\right\}\\[15pt]
\overline{e}[\omega^N,y]&=&\max\limits_x\left\{e^Tx:\begin{array}{l}\,x\in\cX,\\\underline\psi(a[\omega^N]_k,N;y)\leq (A[\omega^N]x)_k\leq\overline\psi(a[\omega^N]_k,N;y),\;k=1,...,\kappa\end{array}\right\}
\end{array}
\end{equation}
where $\underline\psi(\cdot)$ and $\overline\psi(\cdot)$ are defined as in \rf{psifun} of the appendix.
These problems clearly are convex, so  $\underline{e}[\omega^N,y]$ and $\overline{e}[\omega^N,y]$ are efficiently computable. Immediately, we have the following
\begin{lemma}\label{obs} Given $y>1$, the probability of the event
\begin{equation}\label{intervals}
\underline{e}[\omega^N,y]\leq e^T\beta\leq \overline{e}[\omega^N,y],\, ~\forall e,
\end{equation}
is at least $1-2\kappa e{\big(y\big[\ln((y-1)N\big)+2\big]+2\big)}e^{-y}$.
\end{lemma}
\noindent Indeed, when events
\[
a[\omega^N]_k-\overline\psi(a[\omega^N]_k,N;y) \leq F_{\omega^N}(\beta)_k\leq a[\omega^N]_k-\underline\psi(a[\omega^N]_k,N;y),\;k=1,...,\kappa,
\]
take place, $\beta$ is a feasible solution to optimization problems in
(\ref{convopt}). Due to Lemma \ref{lem1002}, this implies that (\ref{intervals})
takes place with probability at least $1-2\kappa e\big(y\big[\ln((y-1)N\big)+2\big]+2\big)e^{-y}$.

\subsection{Estimating parameters of multi-state spatio-temporal processes}\label{multistate}

%\subsection{Recovery for multi-state spatio-temporal processes}\label{multistate}
In this section, we consider the multi-state spatio-temporal process in which an event outcome contains additional information about its category \cite{moller2003statistical}. So far, we considered the case
where at every time instant $t$ every location $k$ maybe be either in the state $\omega_{tk}=0$ (``no event''), or
$\omega_{tk}=1$ (``event''). We are now extending the model by allowing the state of a location at a given time instant to take $M\geq 2$ ``nontrivial'' values
on the top of the zero value ``no event.'' In other words, observation
of the multi-state Bernoulli process is categorical --- we can either observe no event or observe one of $M$ possible event outcomes.
\par
We define {\sl $M$-state spatio-temporal process with memory depth $d$} as follows:
\begin{itemize}

\item We observe a random process on time horizon $\{t:-d+1\leq t\leq N\}$, observation at time $t$ being
$$
\omega_t=\{\omega_{tk}\in\{0,1,\ldots,M\},1\leq k\leq K\}.
$$

\item For every $t\geq 1$, the conditional, $\omega^{t-1}=(\omega_{-d+1},\omega_{-d+2}, \ldots ,\omega_{t-1})$ given, distribution of $\omega_{tk}$ is defined as follows.
With every location $k$, we associate an array of (baseline) parameters $\beta_k=\{\beta_k(p),1\leq p\leq M\}$, and with every pair of locations $k,\ell$
and every $s\in\{1, \ldots, d\}$  --- an array of (interaction) parameters $\beta^s_{k\ell}=\{\beta^s_{k\ell}(p,q), 1\leq p\leq M,0\leq q\leq M\}$.  Then induced by $\omega^{t-1}$ probability of $\omega_{tk}$ to be of category $p$, $1\leq p\leq M$, is given by
    \begin{equation}\label{probstar}
    \textcolor{black}{ \Prob_{\omega^{t-1}}\left\{\omega_{tk}=p\right\}}=\beta_k(p)+\sum_{s=1}^d\sum_{\ell=1}^K\beta^s_{k\ell}(p,\omega_{(t-s)\ell}),
    \end{equation}
     and the probability for $\omega_{tk}$ to take value 0 (no event or ``ground event'') is the
    complementary probability
     \[
    \Prob_{\omega^{t-1}}\left\{\omega_{tk}=0\right\}=1-
     \sum_{p=1}^M\left[\beta_k(p)+\sum_{s=1}^d\sum_{\ell=1}^K\beta^s_{k\ell}(p,\omega_{(t-s)\ell})\right].
     \]
     In other words, $\beta^s_{k\ell}(p,q)$ is the contribution of the location $\ell$ in state $q \in\{0, 1, \ldots, M\}$ at time $t-s$ to the probability for the location $k$ to be in state $p\in\{1, \ldots, M\}$ at time $t$, and $\beta_k(p)$, $p \in \{1, \ldots, M\}$ is the ``endogenous''  component of the probability of the latter event.
\par
Of course, for this description to make sense, the $\beta$-parameters should guarantee that for every $\omega^{t-1}$, that is, for every
collection $\{\omega_{\tau \ell}\in\{0,1, \ldots, M\}:\tau<t,1\leq\ell\leq K\}$, the prescribed by (\ref{probstar}) probabilities are nonnegative and their sum over $p=1, \ldots, M$ is $\leq 1$. Thus, the $\beta$-parameters should satisfy the system of constraints
\begin{equation}\label{newbeta}
\begin{array}{rcl}
0&\leq&\beta_k(p)+\sum_{s=1}^d\sum_{\ell=1}^K\min\limits_{0\leq q\leq M}\beta^s_{k\ell}(p,q),\;1\leq p\leq M,\;1\leq k\leq K,%(a)
\\
1&\geq&\sum_{p=1}^{M } \beta_k(p)+\sum_{s=1}^d\sum_{\ell=1}^K\max\limits_{0\leq q\leq M}\sum_{p=1}^{M}\beta^s_{k\ell}(p,q),\;1\leq k\leq K.
\end{array}
\end{equation}
The solution set $\cB$ of this system is a polyhedral set given by explicit polyhedral representations.
\item We are given convex compact set $\cX$ in the space of parameters $\beta=\{\beta_k,\beta^s_{k\ell}(p,q), 1\leq s\leq d,1\leq k,\ell\leq K,1\leq p\leq M,0\leq q\leq M\}$ such that $\cX$ contains the true parameter $\beta$ of the process we are observing, and $\cX$ is contained in the polytope $\cB$ given by constraints (\ref{newbeta}).
\end{itemize}
We arrange the collection of $\beta$-parameters associated with a $M$-state spatio-temporal process with memory depth $d$ into a column vector (still denoted $\beta$) and denote by $\kappa$ the dimension of $\beta$.\footnote{\label{foot6}In general, $\kappa=KM+dK^2M^2$. However, depending on application, it could make sense to postulate that some of the components of $\beta$ are zeros, thus reducing the actual dimension of $\beta$; for example, we could assume that $\beta_{k\ell}(\cdot,\cdot)=0$ for some ``definitely non-interacting'' pairs $k,\ell$ of locations.}
Note that (\ref{probstar}) says that the $M$-dimensional vector of conditional probabilities for $\omega_{tk}$ to take values $p\in\{1, \ldots, M\}$ given $\omega^{t-1}$ is
$$
[\eta_k^T(\omega_{t-d}^{t-1})\beta]_p
$$
with known to us function $\eta_k(\cdot)$ defined on the set of arrays $\omega_{t-d}^{t-1} \in \{0, 1, \ldots, M\}^{d\times K}$ and taking values in the space of $\kappa\times M$ matrices.
Note that the value of $\omega_{tk}$ is the index of the category, and does not mean magnitude. Same as above, $\eta_k(\omega_{d-1}^{t-1})$ is a Boolean matrix.

To proceed, for $0\leq q\leq M$, let $\chi_q\in\bR^{M}$ be defined as follows: $\chi_0=0\in \bR^{M}$, and $\chi_q$, $1\leq q\leq M$, is the $q$-th vector of the
 standard basis in $\bR^{M}$. In particular, the state $\omega_{tk}$ can be encoded by vector $\bar{\omega}_{tk}=\chi_{\omega_{tk}}$, and the state of our process at time $t$ ---
 by the block vector $\overline{\omega}_t\in\bR^{MK}$ with blocks $\bar{\omega}_{tk}\in\bR^{M}$, $k=1,...,K$.
 In other words: the $k$-th block in $\overline{\omega}_t$ is an $M$-dimensional vector which is the $p$-th basic orth of $\bR^M$ when $\omega_{tk}=p\geq1$, and is the zero vector when $\omega_{tk}=0$. Arranging $\kappa \times M$ matrices $\eta_k(\cdot)$ into a
matrix
$$
\eta(\cdot)=\left[\eta_1(\cdot),...,\eta_K(\cdot)\right] \in \{0, 1\}^{\kappa\times MK},
$$
we obtain
$$
\bE_{|\omega^{t-1}}\left\{\overline{\omega}_t\right\}=\etatt \beta \in \mathbb R^{MK},
$$
where $\bE_{|\omega^{t-1}}$ is the conditional expectation given $\omega^{t-1}$. Note that similarly to Section \ref{sect:prb}, (\ref{probstar}) says that every particular entry in
$\beta$, $\beta_k(p)$ or $\beta^s_{k
\ell}(p,q)$, affects at most one of the entries in the block vector $[\eta_1^T(\omega_{t-d}^{t-1})\beta;...;\eta_K^T(\omega_{t-d}^{t-1})\beta]$ specifically, the $p$-th entry of the $k$-th block, so that the Boolean matrix $\etat $  has at most one nonzero entry in every row.

Note that the spatio-temporal Bernoulli process with memory depth $d$, as defined in Section \ref{sect:prb}, is a special case of $M$-state ($M=1$) spatio-temporal process with memory depth $d$, the case where state 0 at a location contributes nothing to probability of state 1 in another location at a later time, that is, $\beta^s_{k\ell}(1,0)=0$ for all $s,k,\ell$.
\vspace{.1in}
\par\noindent{\em Motivating example: Different types of crime events.} As an illustration, consider a spatio-temporal model of crime events of different types, e.g., burglary and robbery, in a geographic area of interest. We split the area into $K$ non-overlapping cells, which will be our locations. Selecting the time step in such a way that we can ignore the chances for two
or more crime events to occur in the same spatio-temporal cell, we can model the history of crime events in the area as a $M=2$-state spatio-temporal process,
with additional to (\ref{newbeta}) convex restrictions on the vector of parameters $\beta$ expressing our {\it a priori} information on the probability $\beta_{k}(p)$
of a ``newborn'' crime event of category $p$ to occur at time instant $t$
at location $k$ and on the contribution $\beta^s_{k\ell}(p, q)$ of a crime event of category $q$ in spatio-temporal cell $\{t-s,\ell\}$
to the probability of crime event of category $p$, $p\geq1$, to happen in the spatio-temporal cell $\{t,k\}$.

\vspace{.1in}

The problem of estimating parameters $\beta$ of the $M$-state spatio-temporal process from observations of this process can be processed exactly as in the case of the single state spatio-temporal Bernoulli process. Specifically, observations $\omega^N$ give rise to two monotone and affine vector fields on $\cX$, the first observable and the second unobservable:
\begin{equation}\label{eqtwofields}
\begin{array}{rcl}
F_{\omega^N}(x)&=&\underbrace{\left[{1\over N}\sum_{t=1}^N\etat \etatt \right]}_{A[\omega^N]}x
-\underbrace{{1\over N}\sum_{t=1}^N\etat \overline{\omega}_t}_{a[\omega^N]},\\
\overline{F}_{\omega^N}(x)&=&A[\omega^N]x-A[\omega^N]\beta.\\
\end{array}
\end{equation}
The two fields differ only in constant term, $\beta$ is a root of the second field, and the difference of constant terms, same as the vector  $F_{\omega^N}(\beta)$ due to $\overline{F}_{\omega^N}(\beta)=0$, are zero-mean satisfying, for exactly the same reasons as in Section \ref{sec:bound}, concentration bounds \rf{probbnd} and \rf{mart2010} of Lemmas \ref{lem1} and \ref{lem1002}. To recover $\beta$ from observations, we may use the Least Squares (LS) estimate
obtained by solving variational inequality $\VI[F_{\omega^N},\cX]$ with the just defined $F_{\omega^N}$, or, which is the same, by solving
\begin{equation}\label{eq:LS_multiple}
\min_{x\in\cX}\left\{\Psi_{\omega^N}(x):={1\over 2N}\sum_{t=1}^N \|\etatt x-\overline{\omega}_t\|_2^2\right\}.
\end{equation}
Note that (\ref{probbnd}) and \rf{mart2010}, by the same argument as in Section \ref{sec:bound},
imply the validity in our present situation of Theorem \ref{obsini} and Lemma \ref{obs}.

\subsection{Nonlinear link function}\label{nonllink}

So far, our discussion focused on ``linear'' link functions, where past events contribute additively to the probability of a specific event in a given spatio-temporal cell.
We now consider the case of non-linear link functions. This generalizes our model to allow more complex spatio-temporal interactions.

\subsubsection{Single-state process}\label{sec:ssnl}

Let $\phi(\cdot):D\to\bR^K$ be a continuous {\it monotone} vector field defined on a closed convex domain $D\subset \bR^K$ such that
$$
y\in D\Rightarrow 0\leq\phi(y)\leq [1; \ldots; 1].
$$
For example, we may consider ``sigmoid field'' $\phi(u)=[\phi_1(u);...;\phi_K(u)]$ with
$$
[\phi(u)]_k={\exp\{u_k\}\over 1+\exp\{u_k\}},\quad k\leq K, \ D=\bR^K.
$$
Given positive integer $N$, we define a {\sl spatio-temporal Bernoulli process with memory depth $d$ and link function  $\phi$} as
a random process with realizations  $\{\omega_{tk}\in\{0,1\},k\leq K,-d+1\leq t\leq N\}$ in the same way it was done in Section \ref{sect:prb} with
assumptions of Section \ref{sect:prb} replaced with the following:
\begin{itemize}
%\item vector $\beta
%\in\bR^\kappa$ of process's parameters satisfies the restrictions
%\begin{equation}\label{eq31}
%\etatt \beta\in D, \quad \forall 1\leq t\leq N\\
%\end{equation}
%with given functions $\etat $  taking values in the space of $\kappa\times K$ matrices;
\item we are given a convex compact set $\cX \subset\bR^\kappa$ such that the vector of parameters $\beta$ underlying the observed process belongs to $\cX$ and every $\beta\in\cX$ satisfies
\begin{equation}\label{eq31}
\etatt \beta\in D, \quad \forall 1\leq t\leq N
\end{equation}
with given functions $\etat $  taking values in the space of $\kappa\times K$ matrices;
\item the conditional expectation of $\omega_t\in\{0,1\}^K$ given $\omega^{t-1}$ is $\phi(\etatt \beta)$.
\end{itemize}
Let us set
 \begin{equation}\label{old}
\begin{array}{rcl}
F(x)&=&{1\over N}\bE_{\omega^N}\left\{\sum_{t=1}^N\left[\etat \phi\left(\etatt x\right)-\etat \omega_t\right]\right\}:
\cX\to\bR^\kappa,\\
F_{\omega^N}(x)&=&\underbrace{{1\over N}{\sum}_{t=1}^N\etat \phi\left(\etatt x\right)}_{A_{\omega^N}(x)}
-\underbrace{{1\over N}{\sum}_{t=1}^N\etat \omega_t}_{a[\omega^N]}:
\cX\to\bR^\kappa,\\
\overline{F}_{\omega^N}(x)&=&A_{\omega^N}(x)
-\underbrace{{1\over N}{\sum}_{t=1}^N\etat \phi\left(\etatt \beta\right)}_{\overline{a}[\omega^N]}:
\cX\to\bR^\kappa.\\
\end{array}
\end{equation}
We are now essentially in the situation of Section~\ref{process} (where we considered the special case $\phi(z)\equiv z$ of our present situation).
Specifically, $F(\cdot)$ is a monotone (albeit not affine) vector field on $\cX$, $F(\beta)=0$. The empirical version $F_{\omega^N}(x)$, for every $x\in\cX$, is a monotone on $\cX$ vector field which is an unbiased estimate of $F(x)$.
Besides this, $\overline{F}_{\omega^N}(x)$ is a monotone on $\cX$ vector field, and the true vector of parameters $\beta$ underlying our observations solves
the variational inequality $\VI[\overline{F}_{\omega^N},\cX]$(is a root of $\overline{F}_{\omega^N}$). These observations suggest estimating $\beta$ by weak solution to the variational inequality $\VI[F_{\omega^N},\cX]$.

Note that, same as above, vector fields $F_{\omega^N}$ and $\overline{F}_{\omega^N}$ differ only in the constant terms, and this difference is nothing but $F_{\omega^N}(\beta)$ due to $\overline{F}_{\omega^N}(\beta)=0$; moreover $\xi_t=\etat \omega_t-\etatt \beta$ is a martingale difference.
Though deviation probabilities for $F_{\omega^N}(\beta)$ do not obey the same bound
as in the case of $\phi(z)\equiv z$ (since the matrices $\etat $ now not necessarily
are Boolean with at most one nonzero in a row), the reasoning which led us to (\ref{probbnd}) demonstrates
that the vector $F_{\omega^N}(\beta)$ in our present situation does obey the bound
\begin{equation}\label{probbndnew}
%\Prob_{\omega^N}\left\{\|F_{\omega^N}(\beta)\|_\infty>\gamma\Theta / \sqrt{N}\right\}\leq 2\kappa \exp\{-\gamma^2/2\},\quad\forall \gamma\geq0,
\Prob_{\omega^N}\left\{\|F_{\omega^N}(\beta)\|_\infty\geq \Theta\left[\sqrt{\ln(2\kappa/\epsilon)\over 2N}+{\ln(2\kappa/\epsilon)\over 3N}\right]\right\}\leq \epsilon,
\;\forall \epsilon\in(0,1),
\end{equation}
\noindent where $\Theta$ is the maximum, over all possible $\omega_{d-1}^{t-1}$, of the $\|\cdot\|_1$-norm of rows of $\etat $.
Note that in the situation of this section, our $O(1/\sqrt{N})$ exponential bounds on large deviations of $F_{\omega^N}(\beta)$ from zero,
while being good news, do {\sl not} result in easy-to-compute on-line upper-risk bounds and confidence intervals for linear functions of $\beta$.
Indeed, in order to adjust to our present situation Theorem \ref{obsini}, we need to replace the condition numbers $\theta_p[\cdot]$ with constants
of strong monotonicity of the vector field $F_{\omega^N}(\cdot)$ on $\cX$. On the other hand, to adopt the result of Lemma \ref{obs} in the present setting, we need to replace the quantities $\overline{e}$ and $\underline{e}$,
see (\ref{convopt}), with the maximum (resp., minimum) of the linear form $e^Tx$ over the set $\{x\in\cX: \|F_{\omega^N}(x)\|_\infty \leq\delta\}$.
Both these tasks for a {\sl nonlinear} operator $F_{\omega^N}(\cdot)$ seem to be problematic.

\subsubsection{Multi-state processes}\label{nonlinear}

%\paragraph{Multiple states.}
%\label{nonlinear}

The construction in the previous paragraph can be extended to $M$-state processes.
Below, with a slight abuse of notation, we redefine notation for the multi-state processes.

Let us identify two-dimensional $K\times M$ array $\{a_{k\ell}:1\leq k\leq K,1\leq \ell\leq M\}$
with $KM$-dimensional block vector with $K$ blocks $[a_{k1}; a_{k2}; \ldots ;a_{kM}]$, $1\leq k\leq K$,
of dimension $M$ each.
 With this convention, a parametric $K\times M$ array $\psi(z)=\{\psi_{kp}(z)\in\bR:k\leq K,1\leq p\leq M\}$
 depending on $KM$-dimensional vector $z$ of parameters becomes a vector field on $\bR^{KM}$.
 Assume that we are given an array $\phi(\cdot)=\{\phi_{kp}(\cdot)\in\bR:k\leq K,1\leq p\leq M\}$ of the outlined structure such that vector field $\phi(\cdot)$ is continuous and monotone on a closed convex domain $D\subset\bR^{KM}$,
  and for all $y\in D$
\begin{equation}\label{DomainD}
0\leq \phi_{kp}(y)\leq 1,\,1\leq p\leq M,1\leq k\leq K\ \ \&\ \ ~ \sum_{p=1}^{M} \phi_{kp}(y)\leq 1,\quad 1\leq k\leq K.
\end{equation}
We assume that the conditional probability for location $k$ at time $t$ to be in
state $p\in\{1,\ldots,M\}$ (i.e., to have $\omega_{tk}=p$) given $\omega^{t-1}$ is
$$
\phi_{kp}(\etatt \beta)
$$
for some vector of parameters $\beta\in\bR^\kappa$ and known to us function $\eta(\cdot)$
taking values in the space of $\kappa\times KM$ matrices and such that $\eta^T(\omega_{d-1}^{t-1})\beta\in D$
whenever $\omega_{\tau k}\in\{0,1,...,M\}$ for all $\tau$ and $k$. As a result, the conditional probability to have $\omega_{tk}=0$ is
$$
1-\sum_{p=1}^{M}\phi_{kp}(\etatt \beta).
$$
\par
In addition, we assume that we are given a convex compact set $\cX\subset\bR^\kappa$ such that $\beta\in\cX$
and for all such $\beta$
$$
 \etatt \beta \in D,\quad\forall \{\omega_{\tau k}\in\{0,1,...,M\},\,\forall \tau, k\}.
$$
Same as in Section \ref{multistate}, we encode the collection $\{\omega_{tk}:1\leq k\leq K\}$ of locations'
states at time $t$ by block vector $\overline{\omega}_t$ with $K$ blocks of dimension $M$ each, with the $k$-th block
equal to the $\omega_{tk}$-th vector of the standard basis in $\bR^{M}$ when $\omega_{tk}>0$ and equal to 0 when $\omega_{tk}=0$. We clearly have
$$
\bE_{|\omega^{t-1}}\left\{\overline{\omega}_t\right\}=\phi(\etatt \beta).
$$
Setting
\begin{equation}\label{new}
\begin{array}{rcl}
F(x)&=&{1\over N}\bE_{\omega^N}\left\{\sum_{t=1}^N\left[\etat \phi\left(\etatt x\right)-\etat \overline{\omega}_t\right]\right\}:
\cX\to\bR^{\kappa},\\
F_{\omega^N}(x)&=&\underbrace{{1\over N}{\sum}_{t=1}^N\etat \phi\left(\etatt x\right)}_{A_{\omega^N}(x)}
-\underbrace{{1\over N}{\sum}_{t=1}^N\etat \overline{\omega}_t}_{a[\omega^
N]}:
\cX\to\bR^\kappa\\
\overline{F}_{\omega^N}(x)&=&A_{\omega^N}(x)
-\underbrace{{1\over N}{\sum}_{t=1}^N\etat \phi\left(\etatt \beta\right)}_{\overline{a}[\omega^N]}:
\cX\to\bR^\kappa,\\
\end{array}
\end{equation}
(cf. equation (\ref{old})), we can repeat word by word the comment at the end of Section \ref{sec:ssnl}.

\section{Maximum Likelihood estimate}\label{secML estimate}

In the previous sections, we have discussed the Least Squares estimate of the parameter vector $\beta$. Now, we consider commonly used in statistics alternative approach
based on the Maximum Likelihood (ML) estimation. ML estimate is obtained by maximizing
over $\beta\in\cX$ the conditional likelihood of what we have observed, the condition being
the actually observed values of $\omega_{tk}$ for $-d+1\leq t\leq 0$ and $1 \leq k\leq K$. In this section,
we study the properties of the ML estimate and show that its calculation reduces to a convex optimization problem.

\subsection{ML estimation: case of linear link function}\label{secML_setup}

\paragraph{Single state model}
 Assume, in addition to what has been already assumed, that for every $t$ random variables $\omega_{tk}$ are conditionally independent across $k$ given $\omega^{t-1}$. Then the negative log-likelihood, conditioned by the value of $\omega^0$, is given by
{\small$$
L(\beta)={1\over N}\sum\limits_{t=1}^N\sum_{k=1}^K\left[-\omega_{tk}\ln\left(\beta_k+{\sum}_{s=1}^d{\sum}_{\ell=1}^K\beta^s_{k\ell}\omega_{(t-s)\ell}\right)-
(1-\omega_{tk})\ln\left(1-\beta_k-{\sum}_{s=1}^d{\sum}_{\ell=1}^K\beta^s_{k\ell}\omega_{(t-s)\ell}\right)\right].
$$}\noindent
Note that $L(\cdot)$ is a convex function, so the ML estimate in our model reduces to the convex program
\begin{equation}\label{eq:single_ML_estimate}
\min_{x\in \cX} L(x).
\end{equation}
\paragraph{Multi-state model}
Assume that states $\omega_{tk}$ at locations $k$ at time $t$ are conditionally independent across $k\leq K$ given  $\omega^{t-1}$.
Then the ML estimate is given by minimizing, over $\beta\in\cX$,
the conditional negative log-likelihood of collection $\omega^N$ of observations
(the condition being the initial segment $\omega^0$ of the observation).
The objective in this minimization problem is the convex function
\[
L_{\omega^N}(\beta) = -{1\over N}\sum_{t=1}^N\sum_{k=1}^K \psi_{tk}(\beta,\omega^N),
\]
%(this objective is the actual minus log-likelihood divided by $N$),
%\[
%\max\limits_{x\in \cX}\left\{\sum_{t=1}^N\sum_{k=1}^K\psi^{\omega^N}_{tk}(x)\right\},
%\]
where
\begin{equation}
{\psi_{tk}(\beta,\omega^N)}=\left\{\begin{array}{ll} \ln\left([\eta_k^T(\omega_{t-d}^{t-1})\beta]_{\omega_{tk}}\right),&\omega_{tk}\in\{1,\ldots,M\},
\\
\ln\left(1-\sum_{j=1}^{M}[\eta_k^T(\omega_{t-d}^{t-1})\beta]_j\right),&\omega_{tk}=0.
\end{array}\right.
\end{equation}
%$$
%\begin{array}{c}
%\max\limits_{x\in \cX}\left\{\sum_{t=1}^N\sum_{k=1}^K\psi^{\omega^N}_{tk}(x)\right\},\\
%\psi^{\omega^N}_{tk}(x)=\left\{\begin{array}{ll} \ln\left([\eta_k^T(\omega_{t-d}^{t-1})x]_{\omega_{tk}}\right),&\omega_{tk}\in\{1,...,\pi-1\}
%\\
%\ln\left(1-\sum_{\iota=1}^{\pi-1}[\eta_k^T(\omega_{t-d}^{t-1})x]_\iota\right),&\omega_{tk}=0.
%\end{array}\right.\\
%\end{array}
%$$

\paragraph{Toward performance guarantees}
%\noindent
We are about to show that the ML estimate has a structure similar to the LS estimator that we have dealt
within Section \ref{sect1}, and obeys bounds similar to (\ref{probbndnew}).
Given a small positive tolerance $\varrho$, consider $M$-state spatio-temporal process with $K$ locations
and vector of parameters $\beta\in\bR^\kappa$, as defined in Section \ref{multistate}, restricted
to reside in the polyhedral set $B_\varrho$ cut off $\bR^\kappa$ by ``$\varrho$-strengthened''
version of constraints (\ref{newbeta}), specifically, the constraints
\begin{equation}\label{newbetadelta}
\begin{array}{rcl}
\varrho&\leq&\beta_k(p)+\sum_{s=1}^d\sum_{\ell=1}^K\min\limits_{0\leq q\leq M}\beta^s_{k\ell}(p,q),\;1\leq p\leq M\;,1\leq k\leq K,\\
1-\varrho&\geq&\sum_{p=1}^{M -1}\beta_k(p)+\sum_{s=1}^d\sum_{\ell=1}^K\max\limits_{0\leq q\leq M}
\sum_{p=1}^{M}\beta^s_{k\ell}(p,q),\;1\leq k\leq K.\\
\end{array}
\end{equation}
The purpose of strengthening the constraints on $\beta$ is to make the maximum likelihood,
to be defined below, continuously differentiable on the given parameter domain.

In what follows, we treat vectors from $\bR^{KM}$ as block vectors with $K$ blocks of dimension $M$ each.
For such a vector $z$,  $[z]_{kp}$ stands for the $p$-th entry in the $k$-th block of $z$. Let
$$
Z_0=\left\{\omega\in\bR^{MK}: \omega \geq 0,\sum_{p=1}^{M}[\omega]_{kp} \leq 1,\;\forall k\leq K\right\}.
$$
Similarly, for a small positive tolerance $\varrho$, define
\[Z_\varrho=\left\{z\in\bR^{MK}: [z]_{kp}\geq\varrho,\;\forall k,p,\,\sum_{p=1}^{M}[z]_{kp}\leq 1-\varrho,\;\forall k\right\}\subset Z_0.\]
We associate with a vector $w \in Z_0$ the convex function $\cL_w: Z_\varrho\to\bR$,
\begin{equation}\label{cLfun}
\begin{array}{c}
\cL_w(z):=-\sum_{k=1}^K\left[\sum_{p=1}^{M} [w]_{kp}\ln([z]_{kp})+[1-\sum_{p=1}^{M}[w]_{kp}]\ln(1-\sum_{p=1}^{M}[z]_{kp})\right].
\\
\end{array}
\end{equation}
From now on, assume that we are given a convex compact set $\cX\subset B_\varrho$ known to contain
the true vector $\beta$ of parameters. Then the problem of minimizing the negative log-likelihood
 becomes
\begin{equation}\label{minminus}
\min_{x\in\cX}\left\{L_{\omega^N}(x)={1\over N}\sum_{t=1}^N \cL_{\overline{\omega}_t}(\etatt x)\right\},
\end{equation}
where $\overline{\omega}_t=\overline{\omega}_t(\omega^t)$ encodes, as explained in Section \ref{multistate},
the observations at time $t$, and $\etat $ are as defined in Section \ref{multistate}.

Note that by construction, $\overline{\omega}_t$ belongs to $Z_0$. Moreover, by construction,
we have $\eta^T(\omega_{t-d}^{t-1})x\in Z_\varrho$ whenever $x\in B_\varrho$
and $\omega_{tk}\in\{0,1,...,M\}$ for all $t$ and $k$. Now, minimizers of $L_{\omega^N}(x)$ over $x\in \cX$ are exactly the solutions of the variational inequality stemming from $\cX$ and the monotone and smooth vector field (the smoothness property is due to $L_{\omega^N}(x)$ being convex and smooth on $\cX$):
%\begin{equation}\label{FML estimate}
%\begin{array}{c}
\[
F_{\omega^N}(x)=\nabla_x L_{\omega^N}(x)={1\over N}\sum_{t=1}^N\etat
{\theta(\etatt x,\overline{\omega}_t(\omega^t))}
\]
with
\[
\theta(z,\omega)=\nabla_z\cL_{w}(z)=-\sum_{k=1}^K\left[\sum_{p=1}^{M}{[w]_{kp}\over [z]_{kp}}e^{kp}-{1-\sum_{p=1}^{M}[w]_{kp}\over1-\sum_{p=1}^{M}[z]_{kp}}\sum_{p=1}^{M}e^{kp}\right],\quad[w\in Z_0]\\
\]%\end{array}
%\end{equation}
where $e^{kp}\in\bR^{KM}$ is the block-vector with the $p$-th vector of the standard basis in $\bR^{M}$ as the $k$-th block and all other blocks equal to 0.

Note that we clearly have
\begin{equation}\label{notethat}
w\in Z_\varrho\Rightarrow \phi_{w}(w)=0.
\end{equation}
Let us show that
$F_{\omega^N}(\beta)$ is ``typically small'': its magnitude obeys the large deviation bounds
similar to (\ref{probbnd}) and (\ref{probbndnew}). Indeed, let us set $\overline{z}_t(\omega^{t-1})=\etatt \beta$, so that $\overline{z}_t\in Z_\varrho$ due to $\beta\in B_\varrho$.
Invoking (\ref{notethat}) with
$w=\overline{z}_t(\omega^{t-1})$, we have
\[
%\begin{array}{c}
F_{\omega^N}(\beta)={1\over N}\sum_{t=1}^N\underbrace{\etat \vartheta_t[\omega^t]}_{\xi_t},
\]
where
\[\vartheta_t[\omega^t]=-\sum_{k=1}^K\left[\sum_{p=1}^{M}{[\overline{\omega}_t(\omega^t)]_{kp}-[\overline{z}_t(\omega^{t-1})]_{kp}\over [\overline{z}_t(\omega^{t-1})]_{kp}}e^{kp}
+{\sum_{p=1}^{M}\left[[\overline{z}_t]_{kp}-[\overline{\omega}_t(\omega^t)]_{kp}\right]\over 1-\sum_{p=1}^{M}[\overline{z}_t(\omega^{t-1})]_{kp}}\sum_{p=1}^{M}e^{kp}\right].
%\end{array}
\]Since the conditional expectation
of $[\overline{\omega}_t(\omega^t)]_{kp}$ given $\omega^{t-1}$ equals $[\overline{z}_t(\omega^{t-1})]_{kp}$
the conditional expectation of $\xi_t$ given $\omega^{t-1}$ is zero. Besides this, random vectors $\xi_t$
take their values in a bounded set (of size depending on $\varrho$). As a result, $\|F_{\omega^N}(\beta)\|_\infty$
admits bound on probabilities of large deviations of the form (\ref{probbndnew}), with properly selected
(and depending on $\varrho$)  factor $\Theta$. However, for the reasons presented in Section \ref{nonllink}, extracting from this
bound meaningful conclusions on the accuracy of the ML estimate is a difficult task, and it remains an open problem.
%\begin{remark}
\par\noindent{\em Remark} [Decomposition of LS and ML estimation]. In the models we have considered,
the optimization problems  \rf{opt1}, \rf{eq:LS_multiple}, \rf{eq:single_ML_estimate}, and \rf{minminus},  we aim to solve when building the LS and the ML estimates under mild assumptions
are decomposable (in spite of the fact that the observations are  dependent).
Indeed,  vector  \[
\beta=\{\beta_{kp},\beta^s_{k\ell}(p,q),\,1\leq k,\ell\leq K,\,1\leq p\leq M,\,0\leq q\leq M,\,
1\leq s\leq d\}
\] of the model parameters can be split into $K$ subvectors \[
\beta^k=\{\beta_{kp},\beta^s_{k\ell}(p,q),\,1\leq\ell\leq K,\,1\leq p\leq M,\,0\leq q\leq M,\,1\leq s\leq d\}, \;\;k=1,...,K.
\] It is immediately seen that the objectives to be minimized in the problems in question are sums of $K$ terms,
with the $k$-th term depending only on $x^k$. As a result, if the domain $\cX$ summarizing our a priori information on $\beta$ is decomposable: $\cX=\{x: x^k\in\cX_k,1\leq k\leq K\}$, the optimization problems
yielding the LS and the ML estimates are collections of $K$ uncoupled convex optimization problems
in variables $x^k$. Moreover, under favorable circumstances optimization problem \rf{eq:LS_multiple} admits even finer decomposition. Namely, splitting $\beta^k$ into subvectors
\[\beta^{kp}=\{\beta_{kp},\beta^s_{k\ell}(p,q),\,1\leq \ell\leq K,\,1\leq s\leq d,\,0\leq q\leq M\},
\] it is easily seen that the objective in (\ref{eq:LS_multiple}) is the sum, over $k\leq K$ and $p\leq M$,
of functions $\Psi^{kp}_{\omega^N}(x^{kp})$. As a result, when
$\cX=\{x:x^{kp}\in\cX_{kp},1\leq k\leq K,1\leq p\leq M\}$, (\ref{eq:LS_multiple}) is a collection of
$KM$ uncoupled convex problems $\min_{x^{kp}\in\cX_{kp}}\Psi^{kp}_{\omega^N}(x^{kp})$.
\par
The outlined decompositions may be used to accelerate the solution process.

%\end{remark}

\subsection{ML estimate: General link functions}\label{sec:MLest}

%\paragraph{Remark.}
Let us now derive ML estimate for the case of nonlinear link function considered in Section \ref{nonlinear}. In this situation, we strengthen constraints (\ref{DomainD}) on $D$ to
\[%\begin{equation}\label{DomainDnew}
y\in D\Rightarrow \varrho\leq\phi_{kp}(y), ~ \sum_{p=1}^{M} \phi_{kp}(y)\leq 1-\varrho, \quad 1\leq k\leq K, \, 1\leq p \leq M,
\]%\end{equation}
with some $\varrho>0$. Assuming that $\omega_{tk}$'s are conditionally
independent across $k$ given $\omega^{t-1}$, computing ML estimate for the general link-function reduces to solving problem (\ref{minminus}) with $\cL_{w}(z): D\to\bR$, $w\in Z_0$, given by
\[%begin{equation}\label{cLfunnew}
\cL_{w}(z)=-\sum_{k=1}^K\left[\sum_{p=1}^{M} [w]_{kp}\ln(\phi_{kp}(z))+[1-\sum_{p=1}^{M}[w]_{kp}]\ln(1-\sum_{p=1}^{M}\phi_{kp}(z))\right].
\]%end{equation}
Assuming $\phi$ continuously differentiable on $D$ and $\cL_w(\cdot)$ convex on $D$,
we can repeat, with straightforward modifications, everything that was said above
(that is, in the special case of $\phi(z)\equiv z$), including exponential bounds on probabilities
of large deviations of $F_{\omega^N}(\beta)$. However, in general, beyond the case of affine $\phi_{kp}(\cdot)$, function $\cL_w(\cdot)$ becomes nonconvex.
This is due to the fact that convexity on $D$ of functions
$$
-\ln(\phi_{kp}(\cdot)),\;\;-\ln\Big(1-\sum_p\phi_{kp}(\cdot)\Big)
$$
 is a rare commodity. Nevertheless, convexity of these functions does take place in the case
logistic link function
$$
\phi_{kp}(z)={\exp\{a_{kp}(z)\}\over \sum_{q=0}^{M}\exp\{a_{kq}(z)\}}
$$
with functions $a_{kq}(z)$, $0\leq q\leq M$ that are affine in $z$.

\section{Numerical experiments}

\subsection{Experiments with simulated data}\label{sec:sim}

This section presents the results of several simulation experiments
illustrating applications of the proposed Bernoulli process models. We compare performances of Least Squares (LS) and Maximum Likelihood (ML) estimates in terms of $\ell_1$, $\ell_2$, and $\ell_\infty$ norms of the error
of parameter vector recovery. We assume that $d$ (or a reasonable upper bound on it) is known in our simulation examples. The bracket percentage inside the table below shows the norm of the error relative to the norm of the corresponding true parameter vector.

%available at \url{https://github.com/HongtengXu/Hawkes-Process-Toolkit}.})

\subsubsection{Single state spatio-temporal processes}\label{sec:single_experiment}

 First, consider a single state setting with the memory depth $d=8$ and the number of locations $K=8$. The true parameter values are selected randomly from the set $\cX_0$ as follows:
\begin{itemize}
\item $\beta_k\geq 0$, $\beta_{kl}^s \geq 0$; and $\beta_k+\sum_{s=1}^d\sum_{\ell=1}^K\beta^s_{k\ell}\leq 1$, $\forall k$;
\item $\beta^s_{k\ell}=0$  when $|k-\ell|>1$ (interactions are local);
\item For every $1\leq k,\ell \leq K$, $\beta^s_{k\ell}$ is a non-increasing convex function of $s$.\footnote{\label{foot7} Here, the convexity of a function $f(s)$ in $s\in G=\{1,\ldots,d\}$ means that the function is the restriction of a convex function on the segment $[1, d]$ onto the grid $G$ or, which is the same, that $f(s-1)-2f(s)+f(s+1)\geq0, s=2,3,\ldots,d-1$. This translates into the constraint $\beta_{k,\ell}^{s-1} - 2\beta_{k,\ell}^s + \beta_{k,\ell}^{s+1} \geq 0, s=2,3,\ldots,d-1, \forall k,\ell$.}
\end{itemize}
Note that we have imposed additional to \eqref{eq2} constraints on $\beta$.

We report the performance of the LS estimate (obtained by solving $\VI[F_{\omega^N},\cX]$)
and the ML estimate (obtained by solving (\ref{eq:single_ML_estimate})). To ensure a fair comparison, we do not introduce any additional constraints on the interaction coefficients in our estimation procedure (meaning that the LS and ML estimates do not have any prior knowledge about $\cX_0$ and their assumed admissible set $\cX$ is much larger than $\cX_0$).
Utilizing the Matlab implementation \cite{xu2017thap} of the EM algorithm, we also compute estimations of parameters of the commonly used model of Hawkes process with exponential temporal kernel (see, e.g., \cite{zhou2013learning}). The latter is equivalent to assuming that $\beta_{k\ell}^s = a_{k\ell} \tau e^{-\tau s}$, $s = 1, 2, \ldots$, where $\tau > 0$ is the decay rate parameter and $a_{k\ell} > 0$ represents the interactions between two locations.

\begin{figure}[h!]
\centering
\includegraphics[width=0.9\textwidth]{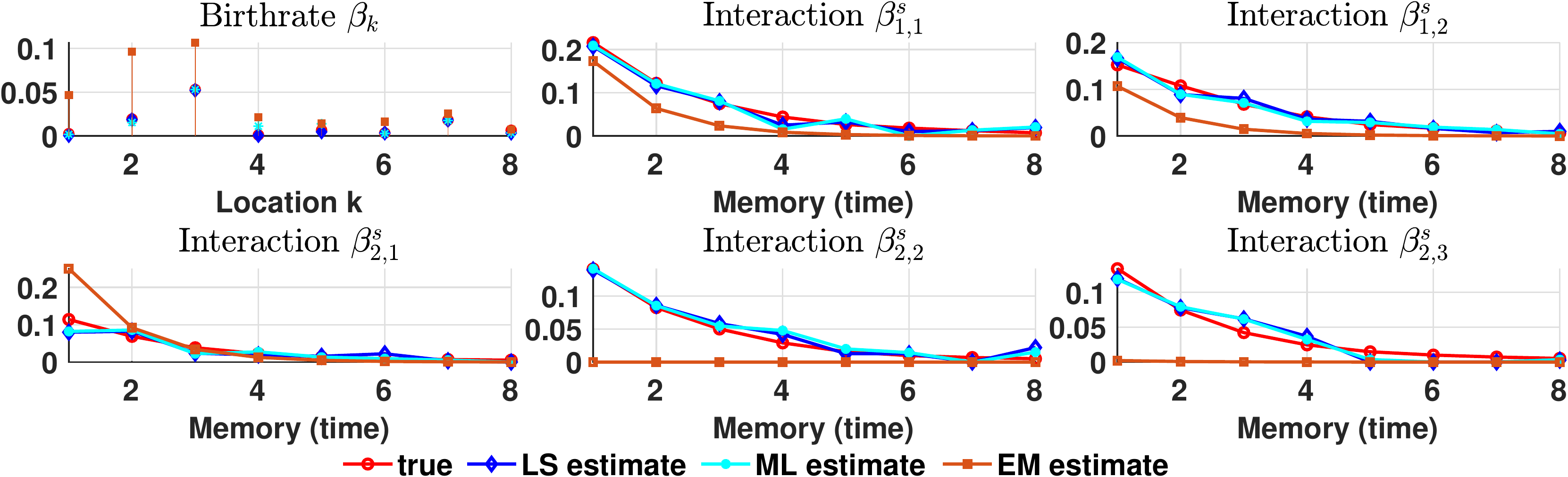} %
\vspace{-0.15in}
\caption{\label{fig:single_recovery} Single-state process: estimates for baseline intensity $\beta_k$ and interactions parameters $\beta_{k\ell}^s$ for one random instance.}
\end{figure}

\begin{figure}[h!]
\centering
\includegraphics[width=0.9\textwidth]{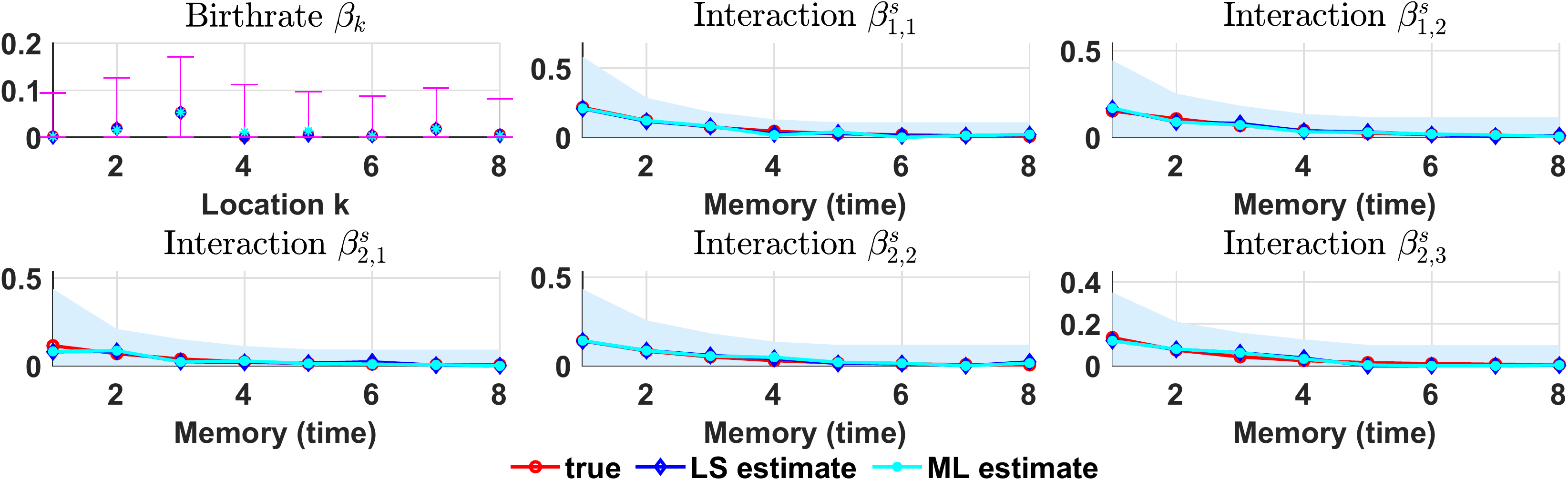}
\vspace{-0.15in}
\caption{\label{fig:single_recovery_1} Computed 90\% confidence intervals corresponding to Fig. \ref{fig:single_recovery}. }
\end{figure}

Fig.\,\ref{fig:single_recovery} shows the recovered interaction coefficients using various methods with $N=10,000$ observations, for a single (randomly generated) instance. The associated error metrics are presented in Table \ref{tab:single_error_one_eg}. The confidence intervals in Fig.\,\ref{fig:single_recovery_1} are computed according to \rf{convopt} by letting $e$ be standard basis vectors in $\bR^\kappa$ and restricting the parameter space to $\cX$. We also repeat the experiment 100 times (each time, generate new true parameters), and the average errors are reported in Table \ref{tab:single_error}. The experiments show that ML and LS estimates exhibit similar performance (ML {outperforming slightly the LS estimates}). Both of them outperform the recovery by EM algorithm based on the exponential kernel, which may be due to a more flexible parameterization of our model.
\begin{table}[h!]\setlength\tabcolsep{12pt}
\centering
\caption{\label{tab:single_error_one_eg} Single-state process:  error of ML, LS, and EM estimation for the one instance shown in Fig. \ref{fig:single_recovery}.}
\renewcommand\arraystretch{1.1}
\begin{tabular}{cccc}
\specialrule{.08em}{0em}{0em}
\hbox{Estimate}& $\ell_1$ error & $\ell_2$ error & $\ell_\infty$ error \\
\hline
%::  below is for original results
%\hbox{ML}& 1.0472 (13.78\%) & 0.1098 (12.91\%)  &0.0375 (12.20\%) \\
%\hline
%\hbox{LS}& 1.1831 (15.57\%) & 0.1237 (14.55\%) &0.0434 (14.10\%)\\
%::  this is new results where we increase significantly the difference between birthrates at different k
%\hbox{ML}& 1.0359 (13.63\%) & 0.1036 (11.93\%)  &0.0261 (10.43\%) \\
%\hline
%\hbox{LS}& 1.1362 (14.95\%) & 0.1091 (12.56\%) &0.0291 (11.59\%)\\
% below is new results which contains comparison with EM
\hbox{ML}& 1.7150 (22.57\%) & 0.1534 (17.67\%)  &0.0342 (13.64\%) \\
\hbox{LS}& 1.8849 (24.80\%) & 0.1714 (19.73\%) &0.0372 (14.84\%)\\
\hbox{EM (exponential kernel)} & 6.3127 (83.06\%) & 0.6413 (73.83\%) &  0.2105 (83.97\%) \\
\specialrule{.08em}{0em}{0em}
\end{tabular}
\end{table}

\begin{table}[h!]\setlength\tabcolsep{12pt}
\centering
\caption{\label{tab:single_error} Single-state process: error of ML, LS, and EM estimation averaged over 100 trials.}
\vspace{-0.1in}
\renewcommand\arraystretch{1.1}
\begin{tabular}{cccc}
\specialrule{.08em}{0em}{0em}
\hbox{Estimate}& $\ell_1$ error & $\ell_2$ error & $\ell_\infty$ error \\
\hline
%::  below is for original results
%\hbox{ML}& 1.0472 (13.78\%) & 0.1098 (12.91\%)  &0.0375 (12.20\%) \\
%\hline
%\hbox{LS}& 1.1831 (15.57\%) & 0.1237 (14.55\%) &0.0434 (14.10\%)\\
%::  this is new results where we increase significantly the difference between birthrates at different k
%\hbox{ML}& 1.0359 (13.63\%) & 0.1036 (11.93\%)  &0.0261 (10.43\%) \\
%\hline
%\hbox{LS}& 1.1362 (14.95\%) & 0.1091 (12.56\%) &0.0291 (11.59\%)\\
% below is new results which contains comparison with EM
%\hbox{ML}& 1.7150 (22.57\%) & 0.1534 (17.67\%)  &0.0342 (13.64\%) \\
%\hbox{LS}& 1.8849 (24.80\%) & 0.1714 (19.73\%) &0.0372 (14.84\%)\\
%\hbox{EM (exponential kernel)} & 6.3127 (83.06\%) & 0.6413 (73.83\%) &  0.2105 (83.97\%) \\
% below is new results averaged over 100 trials
\hbox{ML}& 1.1482 (15.11\%) & 0.1112 (12.60\%) &  0.0336 (11.87\%)  \\
\hbox{LS}& 1.9776 (26.02\%) & 0.1831 (20.72\%) &  0.0472 (16.62\%)\\
\hbox{EM (exponential kernel)} & 6.4725 (85.16\%) & 0.6695 (75.72\%) &  0.2209 (75.17\%) \\
\specialrule{.08em}{0em}{0em}
\end{tabular}
\vspace{-0.15in}
\end{table}

\subsubsection{Multi-state spatio-temporal processes}\label{sec:multi_eg}

Now consider a multi-state spatio-temporal Bernoulli process with the number of states $M=2$. Here the possible states $p=0$ represents no event, $p=1, 2$ represent the event of category 1 and 2, respectively. We assume memory depth $d=8$ and the number of locations $K=10$. The true parameters are randomly generated from the set  $\cX_0$ specified by (again, we impose additional constraints as in Section~\ref{sec:single_experiment}):
\begin{itemize}
\item $\beta_k(p)\geq 0$, $\beta_{kl}^s(p,q) \geq 0$; $\sum_{p=1}^{M}\beta_k(p)+\sum_{s=1}^d\sum_{\ell=1}^K\max_{0\leq q\leq M}\sum_{p=1}^{M}\beta^s_{k\ell}(p,q)\leq 1, \forall k\leq K$;
\item $\beta^s_{k\ell}(p,q)=0$  when $|k-\ell|>1,\ \forall p,q$ (interactions are local);
\item For every $1\leq k,\ell \leq K$ and $1\leq p \leq M,0\leq q\leq M$, $\beta^s_{k\ell}(p,q)$ is a  non-increasing convex function of $s$.
\end{itemize}
Furthermore, we consider two scenarios, with additional constraints on the parameters
\begin{itemize}
\item Scenario 1: events can only trigger future events of the same category, i.e., $\beta^s_{k\ell}(p,q)\equiv0$, $q\neq p$;
\item Scenario 2: events of category $q=0,\ldots,M$,  only trigger events with category $p\leq q$. This can happen, for example, when modeling earthquakes aftershocks: events are marked using $M$ categories according to their magnitudes: $u_1<\ldots < u_{M}$. Set $u_0=0$ and treat the event ``no earthquake'' as ``earthquake of magnitude 0.'' Then each earthquake can trigger ``aftershocks'' with the same or smaller magnitudes. %density: 0.07069
\end{itemize}

\begin{figure}[h!]
\centering
\includegraphics[width=0.9\textwidth]{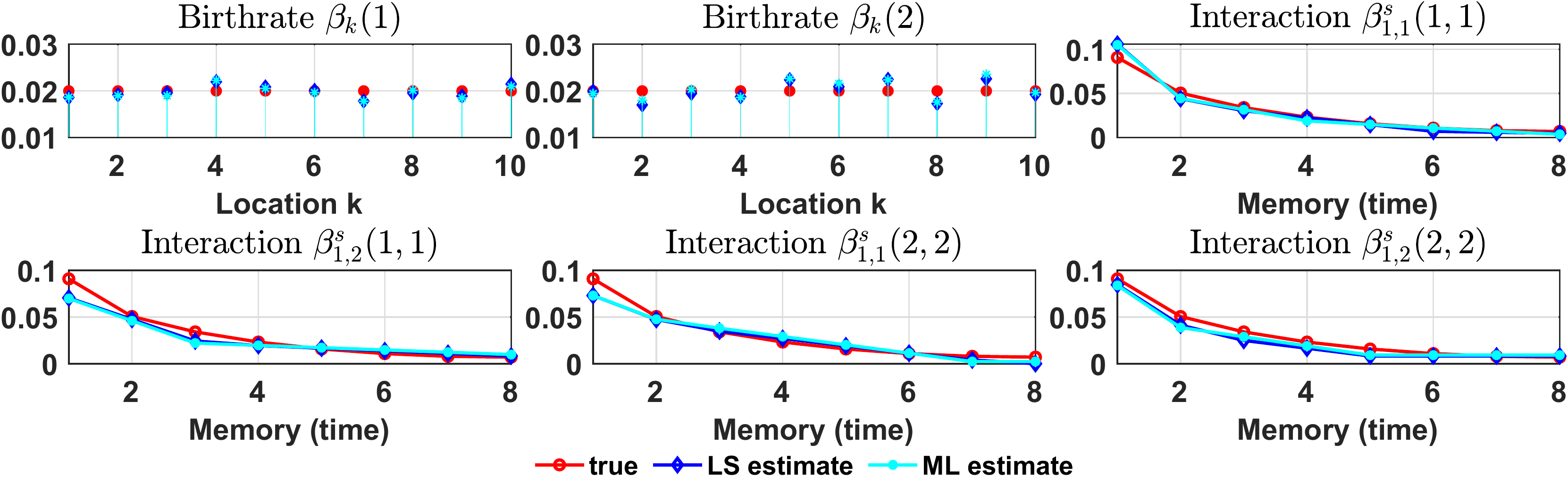}
%\vspace{-0.05in}
\caption{\label{fig:multi_recovery} Multi-state process: examples of LS and ML estimates for baseline intensity $\beta_k(p)$
and interactions parameters $\beta_{k\ell}^s(p,q)$.}
\end{figure}

We generate a synthetic data sequence of length $N=20,000$. For a single (randomly generated) instance, recovery of baseline and interaction parameters are presented  in Fig.\,\ref{fig:multi_recovery}. The associated recovery errors of the LS estimate (solution to \eqref{eq:LS_multiple}) and the ML estimate (solution to  \eqref{minminus}) are reported in Table\,\ref{tab:error_multi}. In addition, we also report the recovery errors separately for (i) the baseline intensity vector (referred to as ``birthrates'')
$\beta_{\text{birth}} = \{\beta_{k}(p),k\leq K, 1\leq p\leq M\}\in\mathbb R^{KM\times 1}$; and (ii) the vector of interactions between different locations
$\beta_{\text{inter}} = \{\beta_{k\ell}^s(p,q)\}\in\mathbb R^{dK^2M(M+1)\times 1}$. As shown in Table \ref{tab:error_multi}, the $\ell_1$ recovery error for estimating birthrate is smaller than that for the interaction parameters. Thus, the recovery error for $\beta$ is dominated by the error for interaction parameters. This could be explained because the magnitude of the baseline intensity is higher than the influence parameters (which is usually needed to have stationary processes).

\begin{table}[h!]\setlength\tabcolsep{12pt}
\centering
\caption{\label{tab:error_multi} Multi-state process recovery: norms of recovery error for LS estimate $\hat\beta_{LS}$ and ML estimate $\hat\beta_{\rm ML}$.
}
\renewcommand\arraystretch{1.1}
\vspace{-0.1in}
\begin{tabular}{ccccc}
\specialrule{.08em}{0em}{0em}
\multirow{2}{*}{Estimate}& \multicolumn{2}{c}{Scenario 1} & \multicolumn{2}{c}{Scenario 2} \\
&$\ell_1$ error  & $\ell_2$ error &$\ell_1$ error &$\ell_2$ error   \\
\hline
$\hat\beta_{\text{ML}}$ & 0.3524 (4.7\%) & 0.0532 (2.5\%) & 1.0179 (13.6\%)  &   0.1146 (5.9\%) \\
$\hat\beta_{\text{LS}}$ & 0.4947 (6.6\%) & 0.0744 (3.4\%) & 1.0854 (14.5\%) &  0.1230 (6.3\%) \\
$\hat\beta_{\text{ML, birth}}$ &  0.0106 (2.7\%) & 0.0028 (3.1\%) & 0.0226 (5.7\%)  & 0.0060 (6.7\%)  \\
$\hat\beta_{\text{LS, birth}}$ &  0.0160 (4.0\%)  & 0.0044 (5.0\%) &0.0237 (5.9\%) &  0.0066 (7.4\%)  \\
$\hat\beta_{\text{ML, inter}}$ & 0.3419 (4.8\%) & 0.0531 (2.5\%) & 0.9952 (14.0\%) &   0.1144 (5.9\%) \\
$\hat\beta_{\text{LS, inter}}$ & 0.4786 (6.7\%) & 0.0743 (3.4\%) & 1.0617 (15.0\%) & 0.1228 (6.3\%)  \\
\specialrule{.08em}{0em}{0em}
\end{tabular}
\end{table}

%\vspace{.2in}

\begin{figure}[h!]
\centering
\begin{tabular}{cc}
\includegraphics[width=0.48\textwidth]{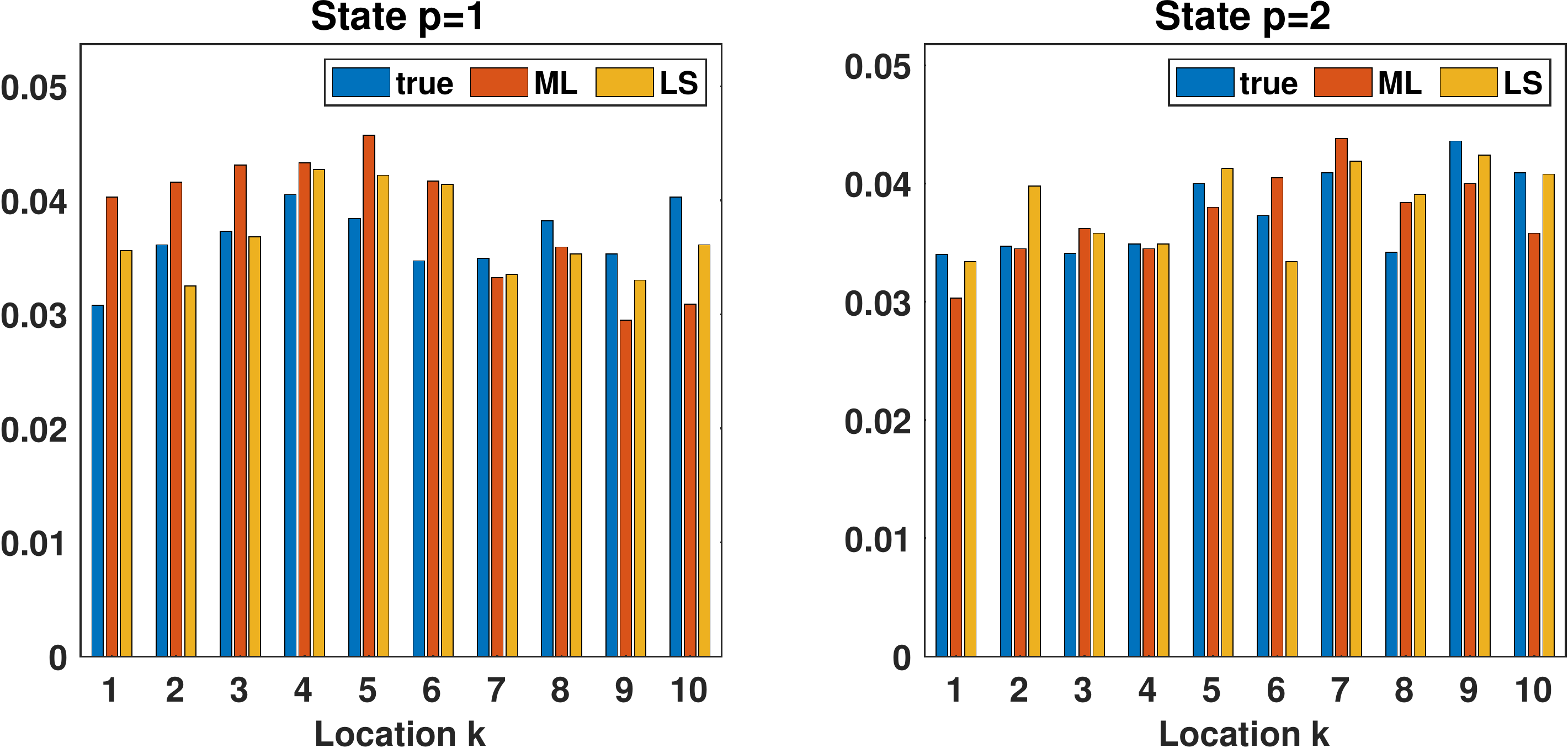} & \includegraphics[width=0.48\textwidth]{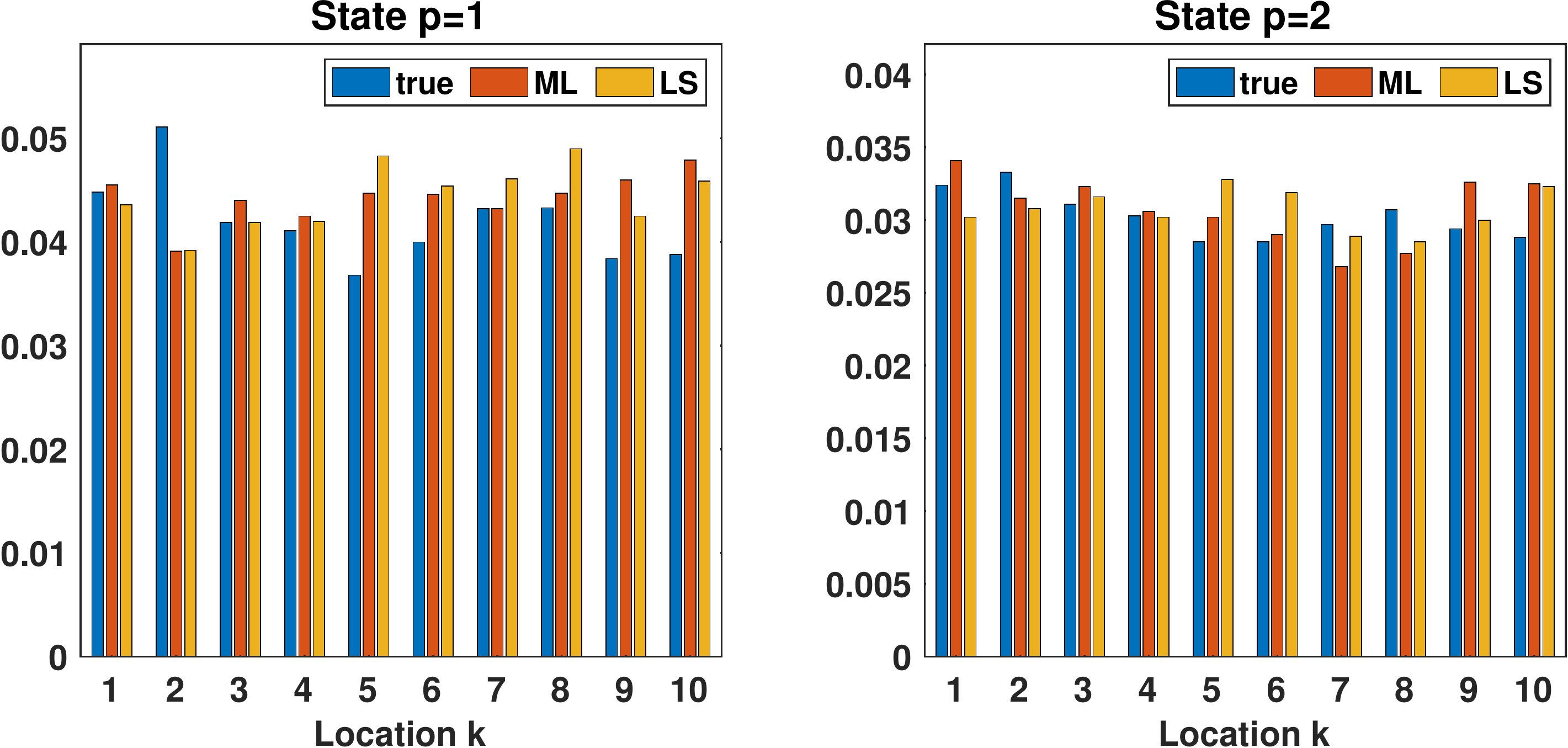}  \\
Scenario 1 & Scenario 2
\end{tabular}
\vspace{-0.15in}
\caption{\label{fre_compare} Multi-state process: experiment to compare the frequency of events from a synthetic sequence (generated using models estimated from training sequence using LS and ML estimates) with that from the testing sequence.}
\vspace{-0.1in}
\end{figure}

Finally, to assess the predictive capability of our model, we did the following experiment. Generate one sequence of discrete events, with length $N = 20, 000$, using randomly selected parameters. We divide the sequence in half: use half for ``training'' and the other half for ``testing''. In particular, we (1) use the first half of the sequence for estimating the Bernoulli process model parameter, (2) use the ``trained'' model to generate a new ``synthetic'' sequence of length $N/2$, and (3) compare the ``synthetic'' sequence with the ``test'' sequence, in terms of the frequency of events, for each category, at each location. The results in Fig.\,\ref{fre_compare} show that the synthetic sequence has a reasonably good match with the testing sequence, based on the LS and the ML estimates.

\subsubsection{Sparse network recovery with negative and non-monotone interactions}

\begin{figure}[h!]
\vspace{-0.05in}
\centering \includegraphics[width=0.35\textwidth]{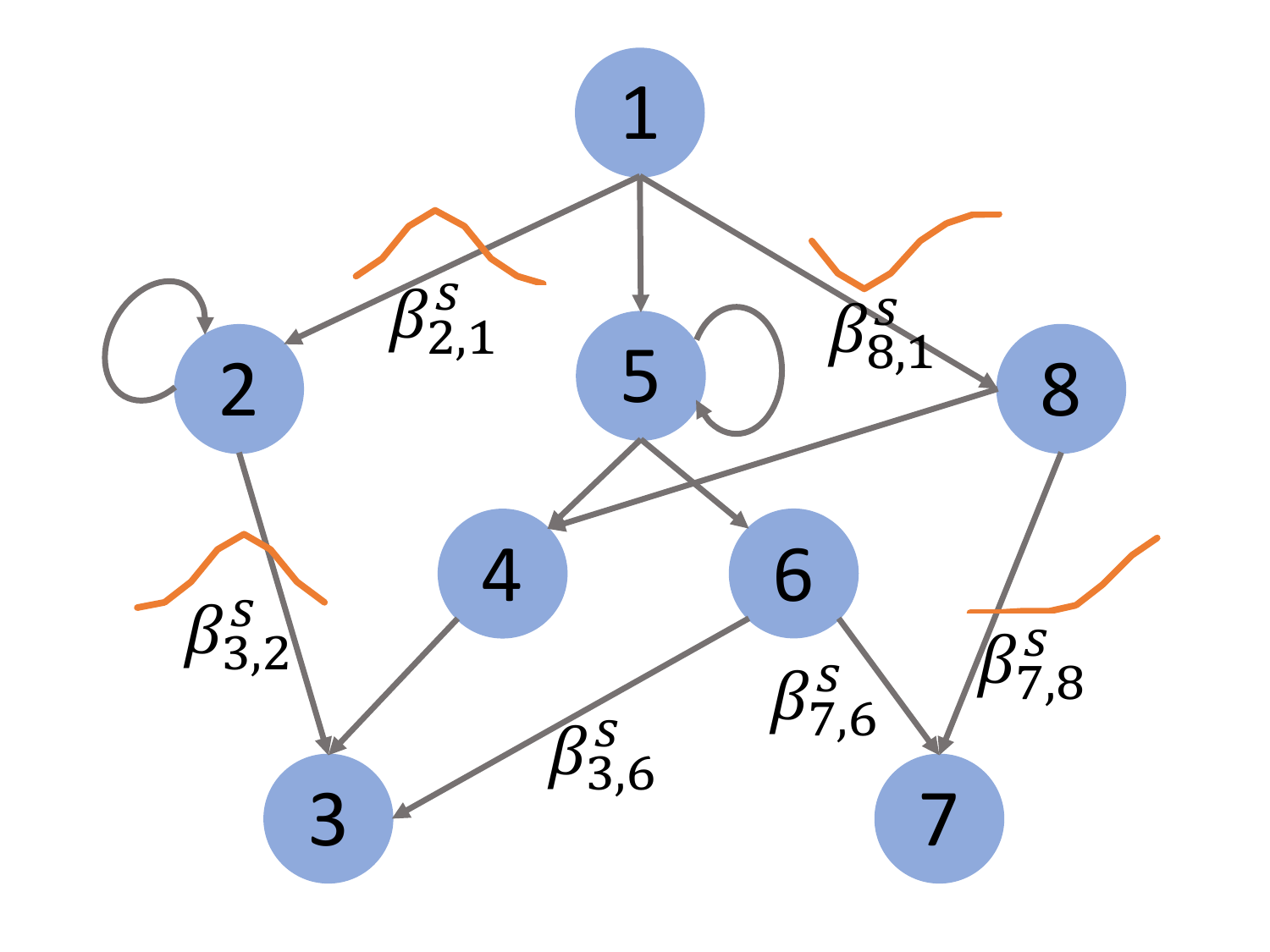}
\vspace{-0.2in}
\caption{\label{fig:graph} Sparse non-planar graph with non-monotonic and negative interaction. Note that the interaction $1\to 8$ is negative.}
\vspace{-0.2in}
\end{figure}

In the last synthetic example, we consider an example to recover a network with ``non-conventional'' interactions: non-monotonic temporal interactions and negative interactions. Consider a sparse, directed, and non-planar graph (meaning that this cannot be embedded on a two-dimensional Euclidean space and, 
thus, this does not correspond to discretized space) with $K = 8$ nodes.
The interaction functions are illustrated in Fig.\,\ref{fig:graph}. 
%\begin{wrapfigure}{r}{0.4\textwidth}
%\end{wrapfigure}

\begin{figure}[h!]
\vspace{-0.15in}
\centering
\includegraphics[width=0.9\textwidth]{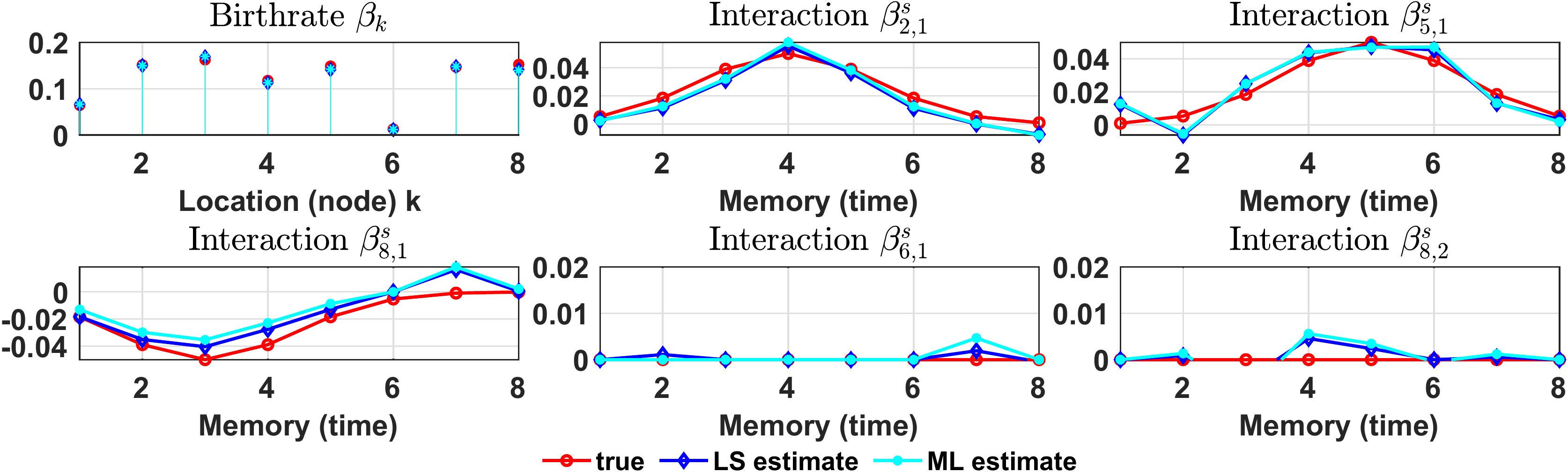}
\vspace{-0.05in}
\caption{\label{fig:graph_recovery} Sparse network identification when graph is unknown: examples of LS and ML estimates of baseline intensity
and vectors of interaction parameters; interactions $\beta_{6,1}$ and $\beta_{8,2}$ correspond to edges $1\to 6$ and $2\to 8$ which do {\it not} exist in the graph in Fig.\, \ref{fig:graph}.}
\vspace{-0.15in}
\end{figure}

\begin{figure}[h!]
\vspace{-0.3in}
\centering
\includegraphics[width=0.35\textwidth]{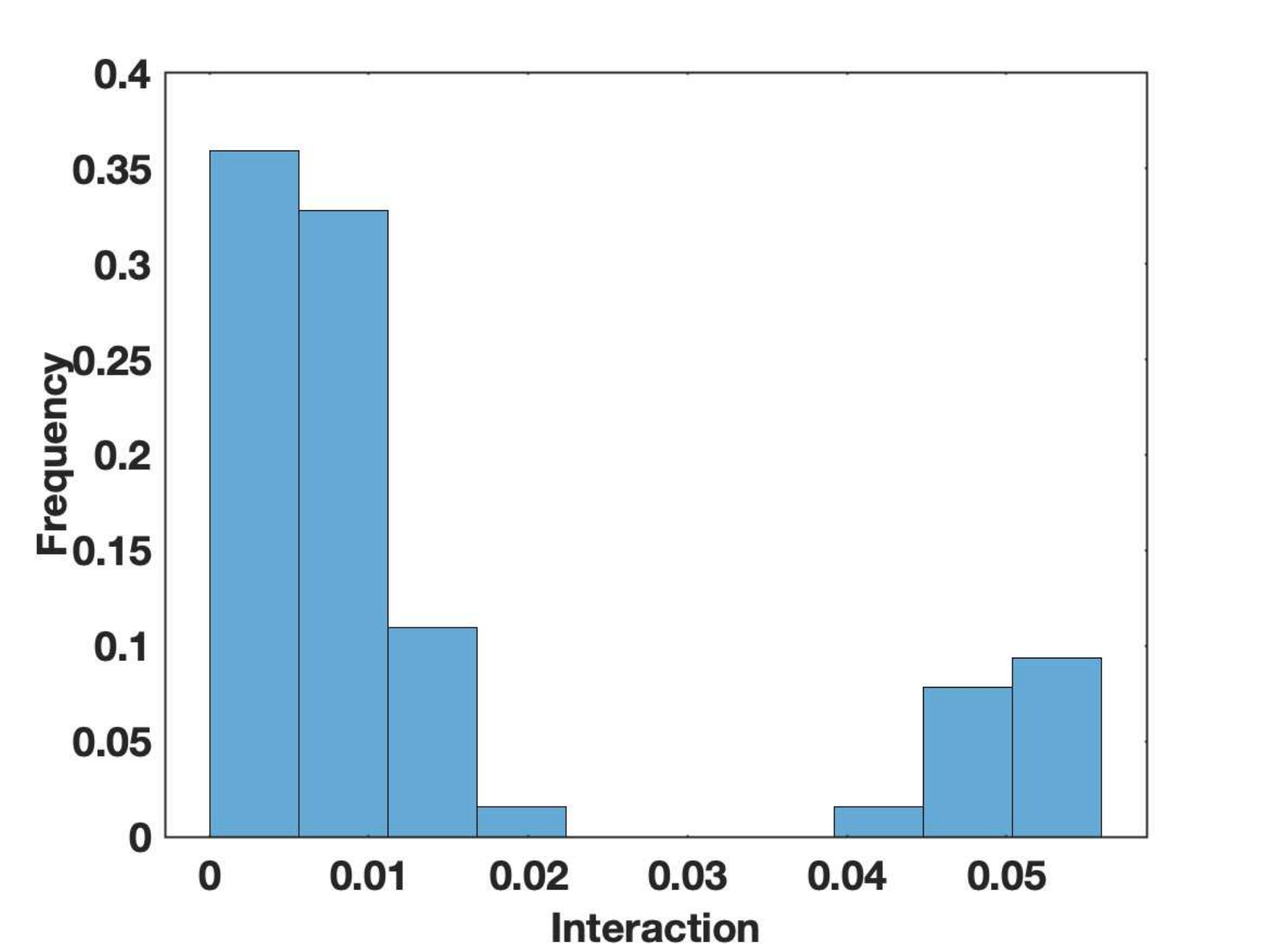}
\vspace{-0.15in}
\caption{\label{fig:graph_hist} Sparse network support recovery: histogram of the recovered interaction parameters $\{ \max_{s=1}^d |\beta_{k,\ell}^s|, 1\leq k,\ell\leq K\}$. Edges with non-zero interactions can be perfectly separated from those with zero interactions. }
%thus, enable us to perform network support recovery without prior knowledge.}
\end{figure}

The baseline intensities are all positive at all 8 nodes. The directed edge (arrows) means there is a one-directional ``influence'' from one node to its neighbor, e.g., $1\to 5$. The self-edges, e.g., $2\to 2$ and $5\to 5$, denote that these nodes have a self-exciting effect: events happen at the node will trigger future events at itself.
The true parameters of the model are generated as follows.
\begin{itemize}
\item Baseline parameters values at
all locations are drawn independently from a uniform distribution on $[0, 0.2]$;
\item For each {\it directed} edge $\ell\rightarrow k$, the interaction $\beta^s_{k,\ell}$ is given by
$\beta^s_{k\ell}=0.05e^{-0.25(s-\tau_{k\ell})^2}$, $s\geq 0$, and the peak $\tau_{k\ell}$ is randomly chosen from $\{1,\ldots,d\}$, except for one edge $1 \to 8$, whose interaction function is set to be negative: $\beta^s_{8,1} = -0.05e^{-0.25(s-\tau_{8,1})^2}$.
\end{itemize}
In our implementation, we consider two scenarios: (1) the graph structure is {\it unknown}: we do {\it not} impose sparsity constraints while obtaining the LS and ML estimates; (2) the graph structure is {\it known}, and then we impose the sparsity constraints by setting the interactions to be 0 when there is no edge; this illustrate the scenario when we have some prior information about the network structure.
We report recovery errors for the two scenarios in Table\,\ref{tab:graph_error} and compare the recovery of interaction parameters under scenario (1) with the true values in Fig.\,\ref{fig:graph_recovery}.
\begin{table}[h!]\setlength\tabcolsep{4pt}
\vspace{-0.1in}
\centering
\caption{\label{tab:graph_error} Sparse network recovery with non-conventional interactions: errors of LS and ML estimates $\hat\beta_{\text{LS}}$, $\hat\beta_{\text{ML}}$.}
\vspace{-0.1in}
\renewcommand\arraystretch{1.1}
%\begin{tabular}{ccccccc}
\begin{tabular}{c@{\hskip 0.2in}ccc@{\hskip 0.25in}ccc}
\specialrule{.08em}{0em}{0em}
\multirow{2}{*}{Estimate}  & \multicolumn{3}{c}{Unknown Graph} &  \multicolumn{3}{c}{Known Graph} \\
& $\ell_1$ error  & $\ell_2$ error & $\ell_\infty$ error & $\ell_1$ error  & $\ell_2$ error & $\ell_\infty$ error \\
\hline
$\hat\beta_{\text{ML}}$ & 1.7694 (58.71\%) & 0.1128 (24.65\%) &  0.0224 (13.79\%)  &
0.4715 (15.64\%) & 0.0593 (12.95\%) &  0.0173 (10.68\%)   \\
$\hat\beta_{\text{LS}}$ &1.8757 (62.23\%) & 0.1166 (25.48\%) &  0.0211 (13.01\%) &
0.4773 (15.84\%) & 0.0606 (13.23\%) &  0.0204 (12.58\%)   \\
$\hat\beta_{\text{ML, birth}}$ &  0.0367 (3.84\%) & 0.0162 (4.42\%) &  0.0111 (6.84\%)    &
0.0126 (1.32\%) & 0.0068 (1.85\%) &  0.0061 (3.75\%)   \\
$\hat\beta_{\text{LS, birth}}$ &  0.0378 (3.95\%) & 0.0172 (4.69\%) &  0.0129 (7.94\%)  &
0.0126 (1.32\%) & 0.0069 (1.89\%) &  0.0061 (3.75\%)   \\
$\hat\beta_{\text{ML, inter}}$ & 1.7327 (84.20\%) & 0.1117 (40.69\%) &  0.0224 (44.73\%)   &
0.4589 (22.30\%) & 0.0589 (21.46\%) &  0.0173 (34.65\%)    \\
$\hat\beta_{\text{LS, inter}}$ &1.8379 (89.31\%) & 0.1153 (42.02\%) &  0.0211 (42.19\%)   &
0.4648 (22.58\%) & 0.0602 (21.92\%) &  0.0204 (40.81\%) \\
\specialrule{.08em}{0em}{0em}
\end{tabular}
\end{table}
From the experiment results, we observe that both the LS and ML estimates match closely with the true parameters, even when the underlying graph structure is unknown. The comparison in Table\,\ref{tab:graph_error} shows a significant improvement in the estimation error when the graph structure is known {\it a priori}. This is consistent with our previous remark that knowing the network structure allows for a better choice of the feasible region resulting in reduced estimation error.
Moreover, by examining the histogram of the maximum interaction between each pair,  i.e., $\{ \max_{s=1}^d |\beta_{k,\ell}^s|, 1\leq k,\ell\leq K\}$ as shown in Fig.\,\ref{fig:graph_hist}, we observe that we can indeed accurately recover the support of the graph: the estimates of the edges with non-zero interactions,
%\begin{wrapfigure}{r}{0.4\textwidth}
%\end{wrapfigure}
are completely separable from the estimates of the edges with zero interactions.
This indicates that we can apply an appropriate threshold (in this case, e.g., $0.03$) to recover precisely the unknown graph structure completely.  This example also shows that even when prior information about the spare structure of the underlying network is not available, LS and ML estimates can recover the underlying network reasonably well, which opens possibilities of applying the proposed approach to perform casual inference \cite{tank2017granger} using discrete-event data.

\subsection{Real data studies: Crime in Atlanta}
\label{sec:rw}

Finally, we study a real crime dataset in Atlanta, USA, to demonstrate the promise of our methods to recover interesting structures from real-data. We consider two categories of crime incidents, ``burglary'' and ``robbery''. These incidents were reported to the Atlanta Police Department from January 1, 2015, to September 19, 2017. The dataset contains 47,245 ``burglary'' and 3,739 ``robbery'' incidents. As mentioned in the introduction, it is believed that crime incidents are related and have ``self-exciting'' patterns: once crime incidence happens, it triggers similar crimes more likely to happen in the neighborhood in the near future \cite{short2008statistical}. Here, we model the data using a multi-state Bernoulli process with two states: no event ($p=0$), burglary ($p=1$), and robbery ($p=2$).

\begin{figure}[h!]
\vspace{-0.1in}
\centering
\begin{tabular}{cc}
\includegraphics[width=0.45\textwidth]{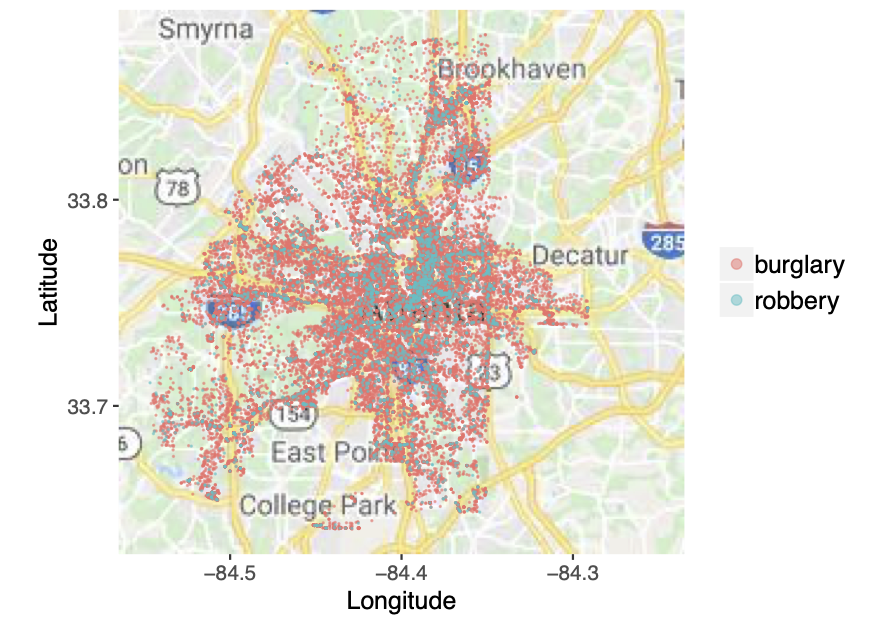} & \includegraphics[width=0.45\textwidth]{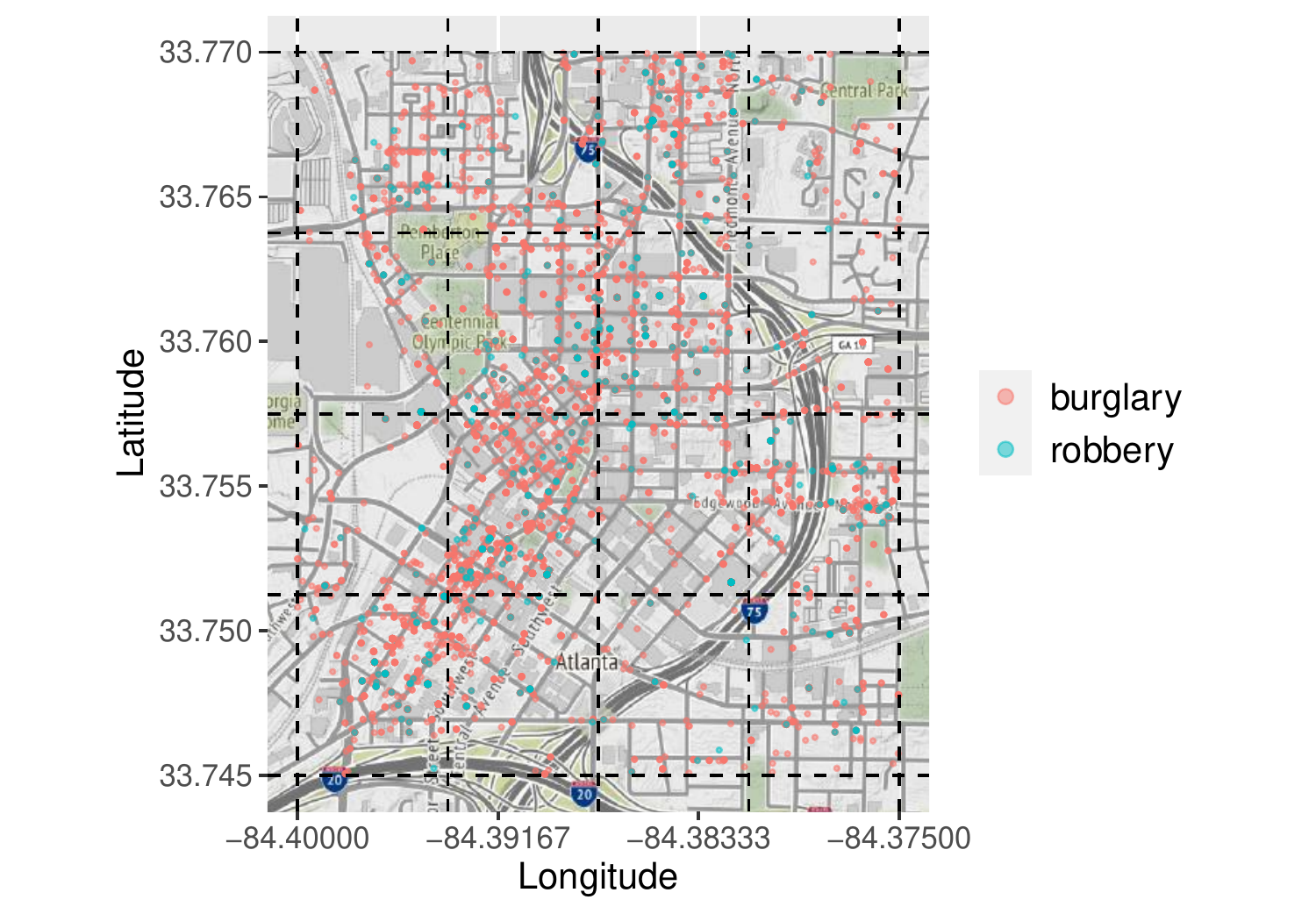}
\end{tabular}
\vspace{-0.1in}
\caption{\label{crime_map} Raw data map: burglary and robbery incidents in Atlanta. Left: the full map; Right: zoom-in around downtown Atlanta.}
\vspace{-0.1in}
\end{figure}
% \fory{~}{Events colors are inverted on the left and the right pictures in Figure \ref{crime_map}.}
We extract crime events around the Atlanta downtown area, as shown in Fig.\,\ref{crime_map}, which contains 6031 ``burglary'' events and 454 ``robbery'' events. The whole time horizon (788 days) is split into discrete time intervals of four hours. The memory depth $d$ is set to $6$ in this example. This value was obtained using a simple ``cross-validation-like'' procedure utilizing predictions of frequencies of the burglary and robbery incidents in various spatial cells. The downtown region is divided uniformly into 16 sub-regions. %; the event timeline is shown in Figure\,\ref{crime_time}.

%\begin{figure}[h!]
%\centering
%\includegraphics[width=0.4\textwidth]{crime_bernoulli.pdf}
%\caption{\label{crime_time} Raw data: times of burglary and robbery events in a 2-year time period. }
%\end{figure}

% \subsection{Recovery and illustration}

\begin{figure}[!ht]
\centering
\begin{tabular}{ccc}
 \multicolumn{3}{c}{Burglary to Burglary} \\
\includegraphics[width=0.25\textwidth]{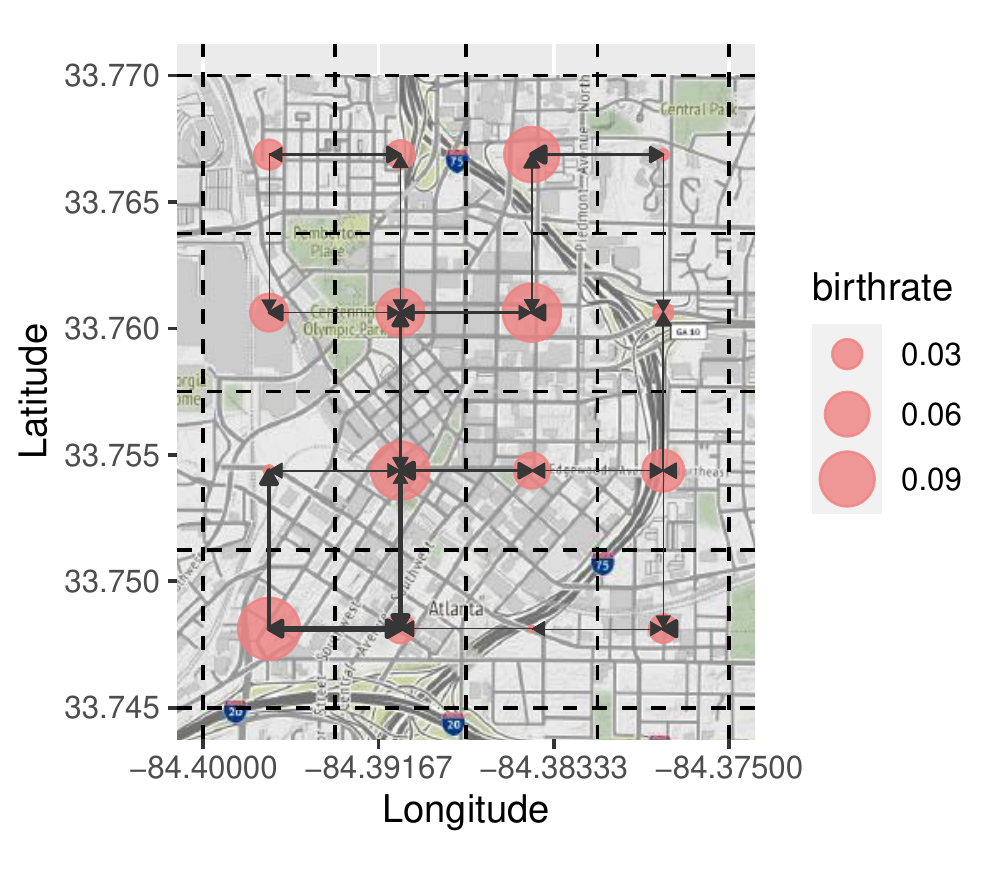} & \includegraphics[width=0.25\textwidth]{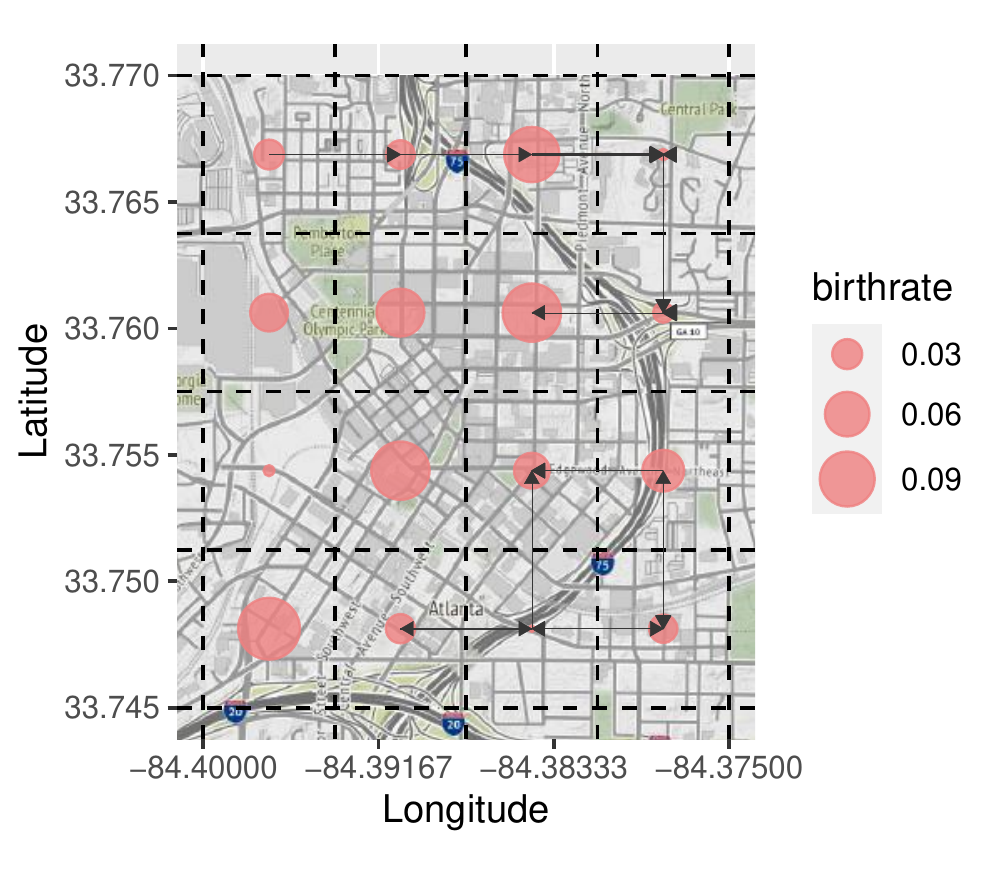} & \includegraphics[width=0.25\textwidth]{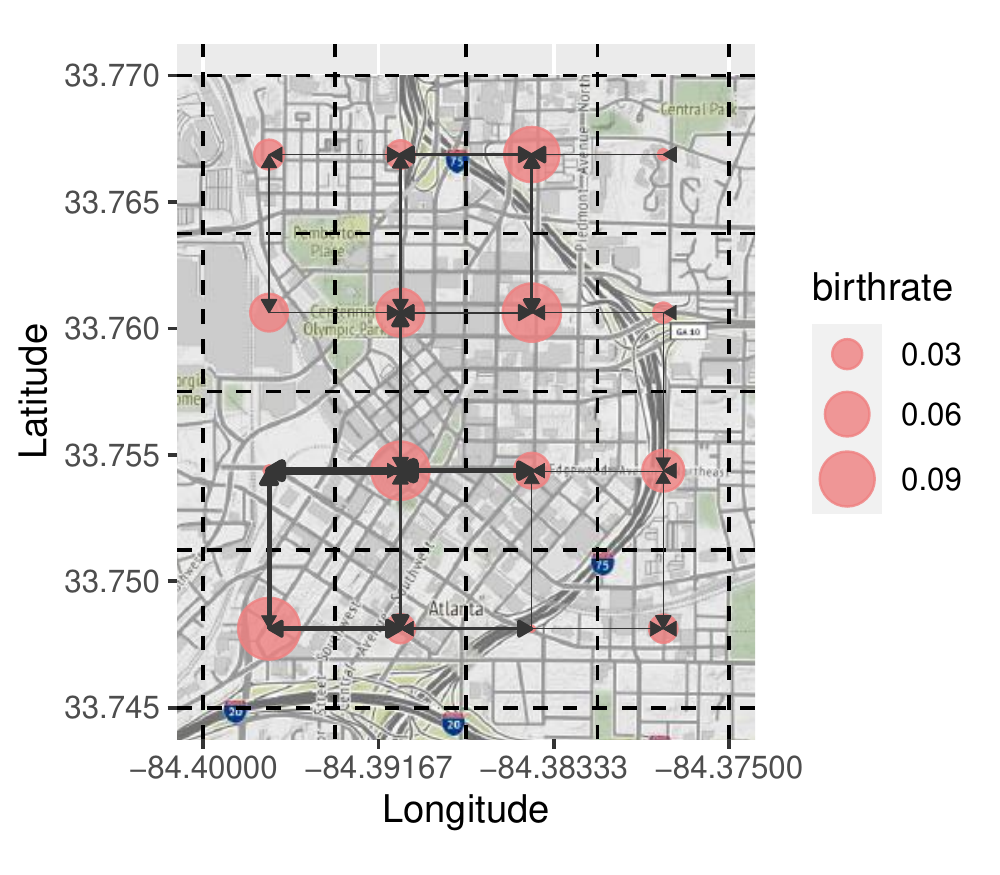} \\
$s=1$ & $s=3$ & $s=6$ \\
%\vspace{0.05in} \\
 \multicolumn{3}{c}{Robbery to Robbery} \\
\includegraphics[width=0.25\textwidth]{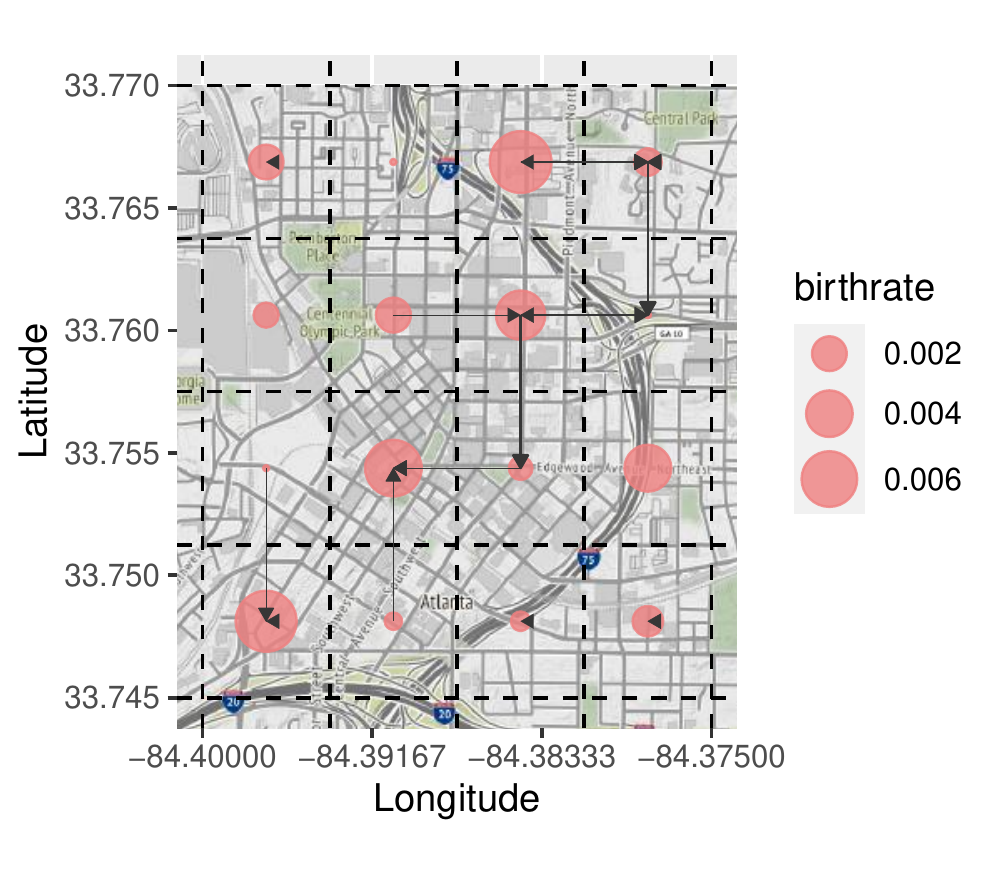} & \includegraphics[width=0.25\textwidth]{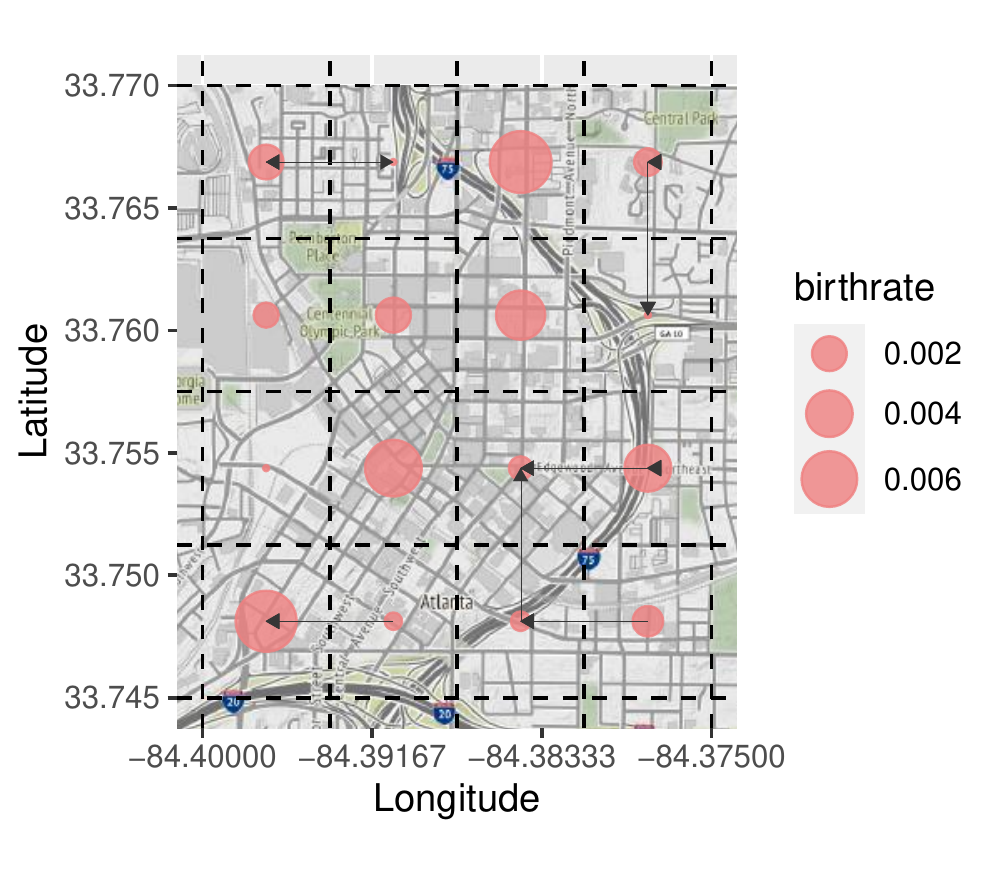} & \includegraphics[width=0.25\textwidth]{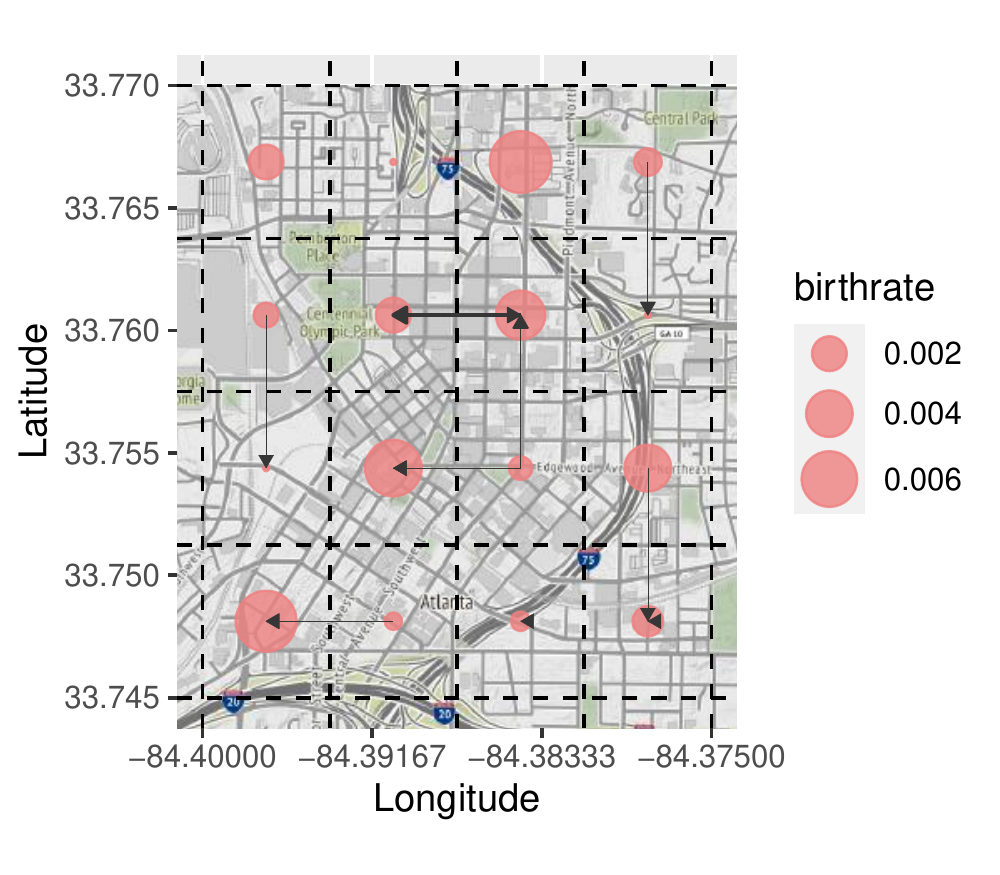} \\
$s=1$ & $s=3$ & $s=6$ \\
 \multicolumn{3}{c}{Burglary to Robbery} \\
\includegraphics[width=0.25\textwidth]{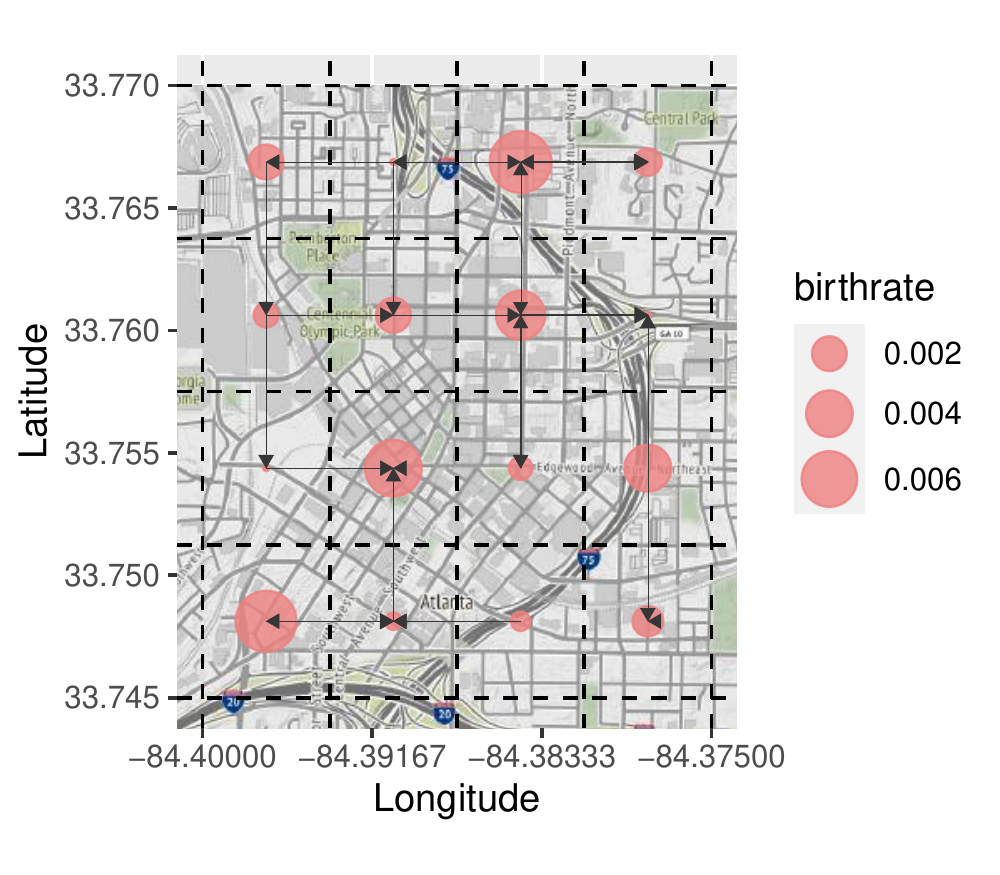} & \includegraphics[width=0.25\textwidth]{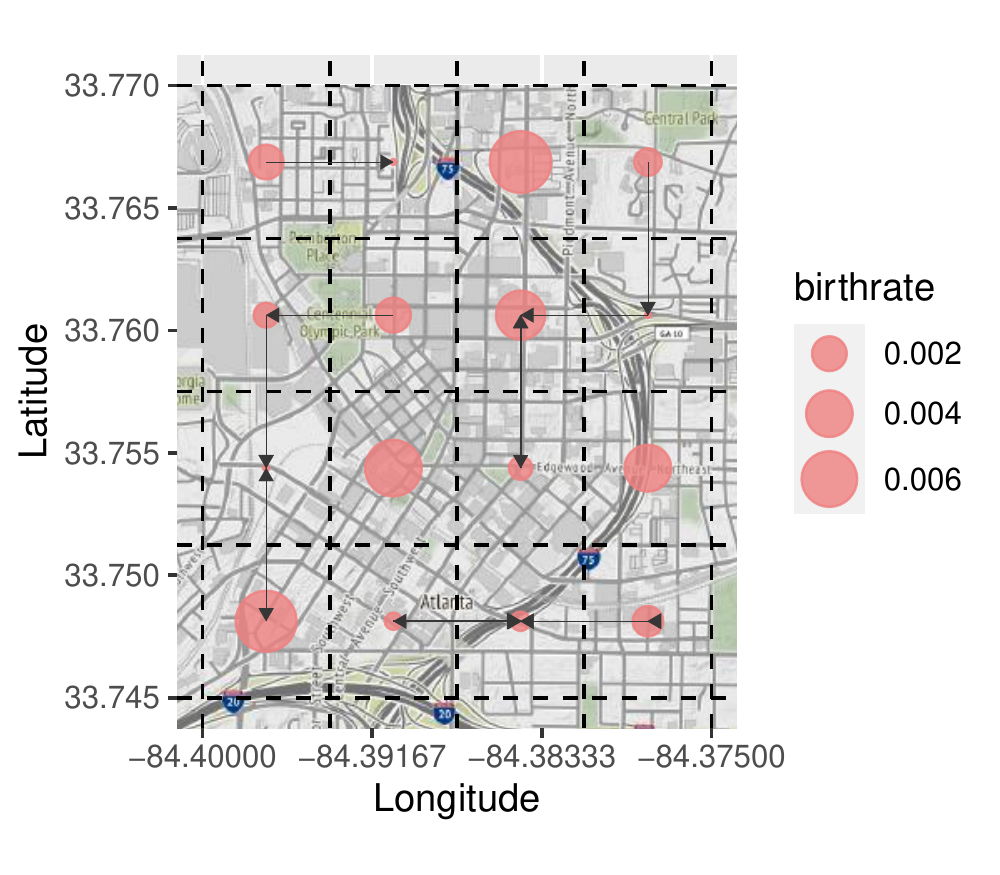} & \includegraphics[width=0.25\textwidth]{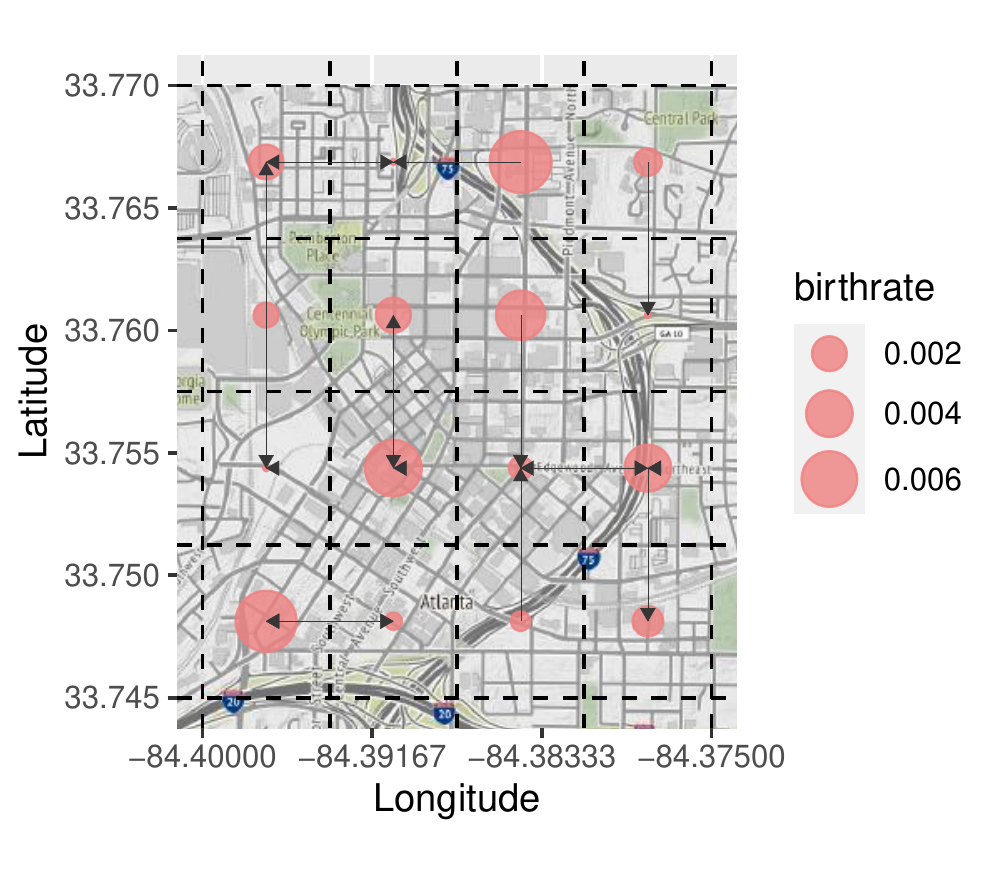} \\
$s=1$ & $s=3$ & $s=6$\\
 \multicolumn{3}{c}{Robbery to Burglary} \\
\includegraphics[width=0.25\textwidth]{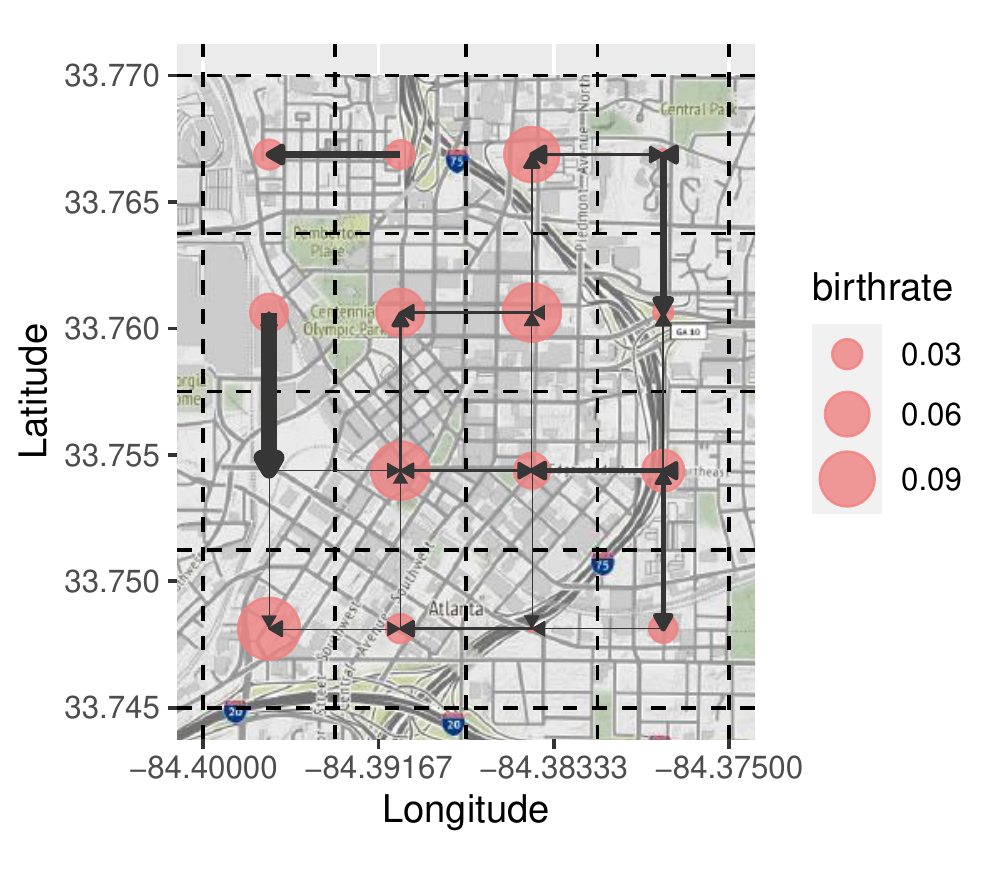} & \includegraphics[width=0.25\textwidth]{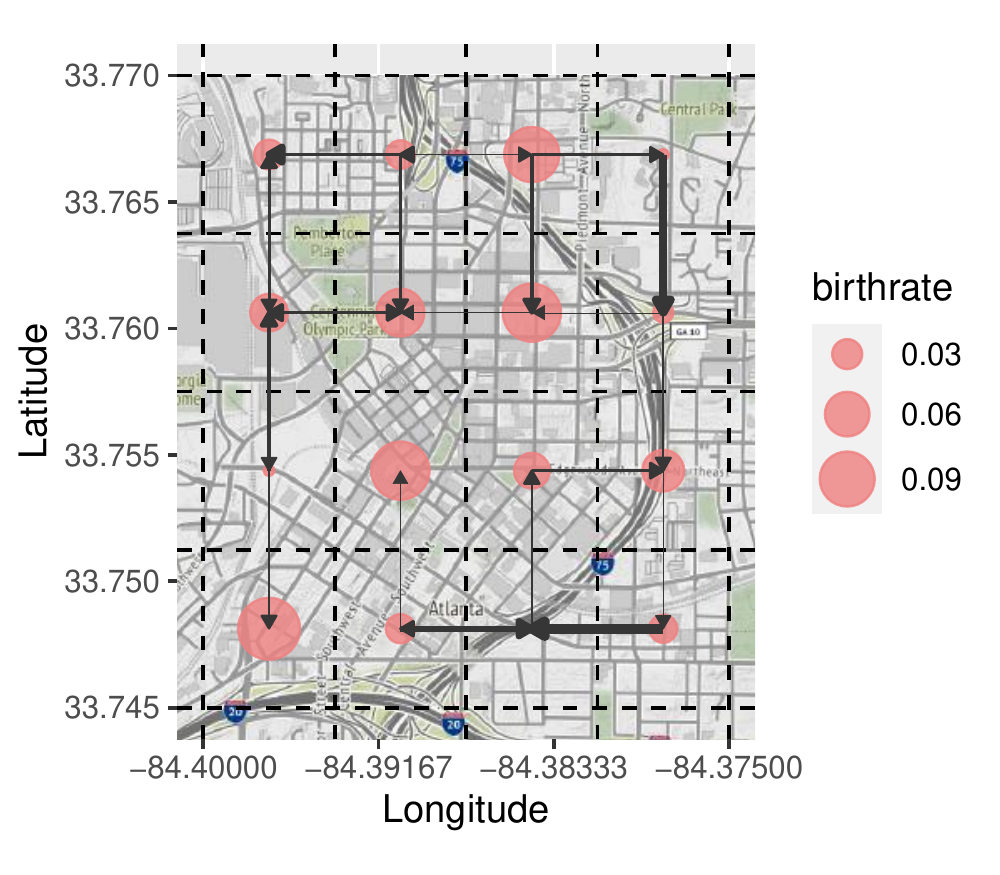} & \includegraphics[width=0.25\textwidth]{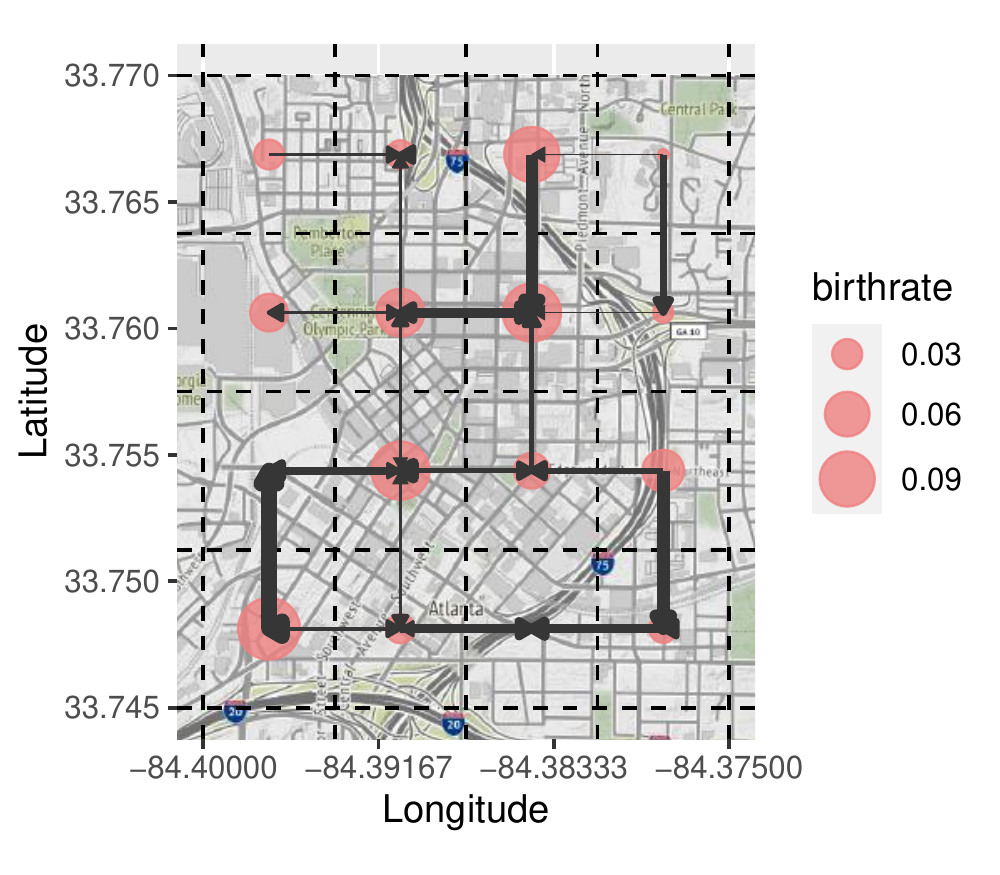} \\
$s=1$ & $s=3$ & $s=6$
\end{tabular}
\caption{\label{crime_aff} Robbery and burglary in downtown Atlanta: recovered spatio-temporal interactions, using LS estimates without additional constraint on the shapes of the interaction functions.}
\end{figure}

We compute the LS estimates of the parameters $\{\beta_k(p),\beta_{k,l}^s(p,q)\}$, in two different ways to set up the constraints: in the first setup, we do not impose additional constraints on the parameters apart from ``basic'' constraints \rf{newbeta}; in the second setup, we impose constraints to only consider temporal interaction function, $\beta^s_{k\ell}$, with monotonic and convex ``shapes''.\footnote{Such constraints are routinely imposed when estimating parameters of
Hawkes model, see, e.g., \cite{Reinhart2017}.} The estimated parameters are shown in Fig.\,\ref{crime_aff}. In the figure, the size of the red dot in each region is proportional to the magnitude of the estimated birthrate $\beta_k(p), k=1,\ldots,K$, for Burglary/Robbery, respectively; the width of the arrow is proportional to the magnitude of the interaction $\beta_{k,l}^s(p,q)$ between locations. It is interesting to notice that our model recovers the dynamic of the interactions and how they change over time. There also seem to be strong interactions between burglary and robbery at different locations.

To validate the model, we experiment similar to we did for the simulated data in Section \ref{sec:multi_eg}. We take the two-year duration of data, divide the sequence in half, use the first half of the sequence to estimate a multi-state Bernoulli process model, generate a synthetic sequence, and compare with the second half of the sequence reserved for testing. We compare the frequency of Burglary and Robbery events across all locations, for the synthetic and testing sequence. The results are shown in Table\,\ref{crime_freq}. The results look to be a reasonably good match, considering that the crime events are relatively rare and with highly complex (and unknown) dynamics: predicting their frequency in the first place is a highly challenging task and an essential research task of criminology.

We also note that the prediction for burglary seems to be better since the frequencies from the synthetic sequence are very close, and the relative error is smaller. This is expected since the number of burglary cases is much larger than the number of robbery cases in our dataset, and the frequency of robbery cases is very small (typically below $0.01$, as shown in Table\,\ref{crime_freq}). The experiment serves as a sanity check and shows that for challenging and noisy real-world datasets, there could be a certain truth to the proposed methods.
\section*{Acknowledgments} Research of Anatoli Juditsky and Arkadi Nemirovski is supported by MIAI \@ Grenoble Alpes (ANR-19-P3IA-0003), CNPq grant 401371/2014-0 and  NSF grant CCF-1523768.
Research of Liyan Xie and Yao Xie are supported by NSF CAREER CCF-1650913, DMS-1938106, DMS-1830210.
\begin{table}[h!]\setlength\tabcolsep{5pt}
\centering
\caption{\label{crime_freq} Crime event model recovery: frequency of Burglary and Robbery events at each location. }
\vspace{-0.1in}
\renewcommand\arraystretch{1.1}
\begin{tabular}{ccccccc}
\specialrule{.08em}{0em}{0em}
\multirow{2}{*}{Locations}  & \multicolumn{3}{c}{Burglary} &  \multicolumn{3}{c}{Robbery } \\
& True & With constr & Without constr & True & With constr & Without constr  \\
\hline
1& 0.1499 &\bf{0.1707} & 0.1766  & 0.0102 & 0.0195  & \bf{0.0186}  \\
2& 0.0284& \bf{0.0373}&0.0445  & 0.0017& \bf{0.0203}  &0.0212 \\
3 & 0.0483& \bf{0.0580} &0.0606  & 0.0021& 0.0254  &\bf{0.0195}  \\
4 & 0.0407& \bf{0.0364} &0.0356  & 0.0017& \bf{0.0178} &0.0224  \\
5 & 0.0508& \bf{0.0529} &0.0648  & 0.0042&0.0220 &\bf{0.0165} \\
6 & 0.1957& 0.2088 &\bf{0.1834}  & 0.0131&0.0208 &\bf{0.0144} \\
7 & 0.0970&0.1368 &\bf{0.1224}  & 0.0068&0.0229 &\bf{0.0191} \\
8 & 0.0419& 0.0580 &\bf{0.0563}  & 0.0021&\bf{0.0127} &0.0182 \\
9 & 0.0148& \bf{0.0161} &0.0220  & 0.0013&\bf{0.0165} &0.0212 \\
10 & 0.0584& \bf{0.0729} &0.0805  & 0.0055&0.0258 &\bf{0.0178} \\
11 & 0.1266& \bf{0.1525} &0.1529  & 0.0106&0.0195 &\bf{0.0169} \\
12 & 0.1364 & \bf{0.1266} &0.1186  & 0.0102&0.0191 &\bf{0.0169} \\
13 & 0.0322& 0.0521 &\bf{0.0445}  & 0.0021&0.0229 &\bf{0.0224} \\
14 & 0.0627& 0.0868 &\bf{0.0834}  & 0.0055&0.0212 &\bf{0.0195} \\
15 & 0.0208& \bf{0.0224} &0.0280  & 0.0008&0.0241 &\bf{0.0216} \\
16 & 0.0144&0.0203 &\bf{0.0178}  & 0.0013&\bf{0.0203} &0.0212 \\
\specialrule{.08em}{0em}{0em}
\end{tabular}
\vspace{-0.2in}
\end{table}

%

%\clearpage
\bibliographystyle{IEEEtran}
\bibliography{ReferencesChP-0430}

\appendix
\small
\section{Proof of Lemma \ref{lem1002}
%(Bernstein's inequality for Bernoulli martingales)
}
We start with describing an application of the Bernstein inequality for martingales (cf., e.g., \cite{azuma1967weighted,freedman1975tail,fan2012hoeffding,bercu2015concentration}) in our situation. Let $\omega_i\;i=...,0,1,2,...$ be a sequence of random binary vectors in $\bR^m$ such that the conditional distribution of the $j$-th component $\omega_{ij}$, $j=1,...,m$, of $\omega_i$ given $\omega^{i-1}$ is Bernoulli distribution with parameter $\mu_{ij}=\bE_{|\omega^{i-1}}\{\omega_{ij}\}$. Now, consider the sequence of Boolean vectors $\gamma_i,\,i=1,2,...$, $\gamma_i\in \bR^m$, such that $\gamma_i$ is $|\omega^{i-1}$-measurable with $\sum_j\gamma_i^j\leq 1$ a.s..
Finally, let $\zeta_i=\gamma_i^T\omega_i-\gamma_i^T\mu_i$; note that, in this case,
\[
\bE_{|
\omega^{i-1}}\{\zeta_i\}=0,\;\sigma_i^2:=\bE_{|\omega^{i-1}}\{\zeta_i^2\}=\gamma_i^T\mu_i(1-\gamma_i^T\mu_i)\leq \four,\;\mbox{and}\;|\zeta_i|\leq 1\;\;\mathrm{a.s..}
\]
%Let $I=\sum_{i=1}^N\sum_{j=1}^m \gamma_{i}^j$, the number of nonzero vectors among $\gamma_1,...,\gamma_N$; we
Denote $\bar \mu_N={1\over N}\sum_{i=1}^N \gamma_i^T\mu_i$, $\bar \nu_N={1\over N}\sum_{i=1}^N \gamma_i^T\omega_i$, $\bar s_N={1\over N }\sum_{i=1}^N \sigma_i^2$, and  $\bar\zeta_N={1\over N }\sum_{i=1}^N \zeta_i$.
%\par
%The following result is an immediate application of the Bernstein inequality for martingales.
\begin{lemma}\label{lem:mart}
Let $0<\underline s<\overline s<\infty$, and let $y> 1$. One has
\be
\Prob\left\{|\bar \zeta_N|\geq \sqrt{2y\bar s_N\over N }+{ y\over 3N },\;\underline s\leq \bar s_N\leq \overline s\right\}\leq 2e\big(y\ln\big({\overline s/ \underline s}\big)+1\big)e^{-y}.
\ee{mart1}
and, as a consequence,
\be
\Prob\left\{|\bar\zeta_N|\geq \sqrt{2y\bar s_N\over N }+{y\over 3N }\right\}\leq 2e \big(y\big[\ln((y-1)N\big)+2\big]+2\big) e^{-y}.
\ee{mart10}
Moreover, we have
\be
\Prob\left\{\underline\psi(\bar \nu_N,N;y) \leq \bar\mu_N\leq \overline\psi(\bar \nu_N,N;y)\right\}\geq 1-2e \big(y\big[\ln((y-1)N\big)+2\big]+2\big)  e^{-y}
\ee{mart0100}
where
\be
\begin{array}{rcl}
\underline\psi(\nu,N;y)&=&\left\{
\begin{array}{l}
(N+2y)^{-1}\left[N\nu+{2y\over 3}-\sqrt{2N\nu y+{y^2\over 3}-{2y\over N}\left({y\over 3}-\nu{N}\right)^2}\right]\;\mbox{if}\;\nu>{y\over 3N},\\
0\quad\mbox{otherwise;}\end{array}\right.\\
\overline\psi(\nu,N;y)&=&\left\{\begin{array}{l}
(N+2y)^{-1}\left[N\nu+{4y\over 3}+\sqrt{2N\nu y+{5y^2\over 3}-{2y\over N}\left({y\over 3}+\nu{N}\right)^2}\right]\;\mbox{if}\;\nu<1-{y\over 3N},\\
1\quad\mbox{otherwise,}\end{array}\right.
\end{array}\ee{psifun}
so that
\be
\Prob\left\{\bar \nu_N-\overline\psi(\bar \nu_N,N;y) \leq \bar\zeta_N\leq \bar\nu_N-\underline\psi(\bar \nu_N,N;y)\right\}\geq 1-2e {\big(y\big[\ln((y-1)N\big)+2\big]+2\big)} e^{-y}.
\ee{mart010}
\end{lemma}
\par\noindent{\em  Proof of the lemma.}
Utilizing Bernstein's inequality for martingales (cf., e.g., \cite[Theorem 3.14]{bercu2015concentration}) we obtain for all $z> 0$ and $s> 0$,
\be
\Prob\left\{\left|\sum_{i=1}^N\zeta_i\right|\geq \sqrt{2zs}+{z\over 3},\,\sum_{i=1}^N\sigma_i^2\leq s,\right\}\leq 2e^{-z}.
\ee{BDR}
We conclude that
\[
\Prob\left\{|\bar \zeta_N|\geq \sqrt{{2\bar s_N\over N }z(1+ z^{-1})}+{z\over 3 N },\, \bar s_N\in[s, (1+z^{-1})s]\right\}\leq 2e^{-z},
\]
implying that for $y=z+1> 1$
\be
\Prob\left\{|\bar\zeta_N|\geq \sqrt{2y\bar s_N\over N }+{y
\over 3N },\,\bar s_N \in\left[s,{( y-1)^{-1}ys}\right]\right\}\leq 2e^{-y+1}.
\ee{bythis0}
Let now $s^j=\min\left\{\overline s,\left(y\over y-1\right)^j s^0\right\}$, $j=0,..., J$, with  $s^0=\underline s$, $s^{J}=\overline s$, and
%\[
$J=\left\rfloor\ln\big({\overline s/ \underline s}\big)\ln^{-1}\left(( y-1)^{-1}y\right)\right\lfloor$.
%\]
Note that $\ln\big(1+1/(y-1)\big)\geq 1/y$ for $y>1$, so that
\[J \leq \ln\big({\overline s/ \underline s}\big)\ln^{-1} \big(( y-1)^{-1}y\big)+1\leq y\ln\big({\overline s/ \underline s}\big)+1.
\]
On the other hand, due to \rf{bythis0},
%\bse
%\lefteqn{\Prob\left\{|\bar \zeta_N|\geq \sqrt{2y\bar s_N\over N }+{y\over 3N },\;\underline s\leq \bar s_N \leq \overline s\right\}}\\&\leq&\sum_{j=1}^{J}
%\Prob\left\{|\bar \zeta_N|\geq \sqrt{2y\bar s_N\over N }+{y
%\over 3N },\,\bar s_N \in\left[s^j,s^{j+1}\right]\right\}\\&
%\leq& 2J e^{-y+1}\leq 2e\big(y\ln\big({\overline s/ \underline s}\big)+1\big)e^{-y}%\qquad\qquad\mbox{\qed}
%\ese
\[
\begin{aligned}
\Prob\left\{|\bar \zeta_N|\geq \sqrt{2y\bar s_N\over N }+{y\over 3N },\;\underline s\leq \bar s_N \leq \overline s\right\}&\leq\sum_{j=1}^{J}\Prob\left\{|\bar \zeta_N|\geq \sqrt{2y\bar s_N\over N }+{y
\over 3N },\,\bar s_N \in\left[s^j,s^{j+1}\right]\right\}\\
&\leq 2J e^{-y+1}\leq 2e\big(y\ln\big({\overline s/ \underline s}\big)+1\big)e^{-y}
\end{aligned}
\]
what is \rf{mart1}.
Let us put $s=(18z)^{-1}$ in \rf{BDR}; together with $y=z+1> 1$, we get
%\[
%\Prob\left\{|\bar \zeta_N|\geq {z+1\over 3N },\,\bar s_N \leq {1\over 18Nz}\right\}\leq 2e^{-z},
%\]
%or
\be
\Prob\left\{|\bar \zeta_N|\geq {y\over 3N },\,\bar s_N \leq {1\over 18N(y-1)}\right\}\leq 2e^{-y+1}.
\ee{mart2}
Furthermore, we have $\bar s_N \leq 1/4$ a.s.. When substituting
$\underline s=(18(y-1))^{-1}$ and $\overline s=N/4$ into \rf{mart1} we obtain
\[
\Prob\left\{|\bar \zeta_N|\geq \sqrt{2y\bar s_N\over N }+{y\over 3N },\,\bar s_N \geq {1\over 18N(y-1)}%\leq {\bar \gamma^2 n\over 4}
\right\}\leq
2e\big(y\ln\big(\mbox{\small ${9\over 2}$}(y-1)N\big)+1\big)e^{-y}.
\]
Finally, when taking into account \rf{mart2} we conclude with
$$
\Prob\left\{|\bar \zeta_N|\geq \sqrt{2y\bar s_N\over N }+{y\over 3N }%,\,{\bar \gamma^2\over 18(y-1)}\leq s_n\leq {\bar \gamma^2 n\over 4}
\right\}\leq 2e\big(y\ln\big(\mbox{\small ${9\over 2}$}(y-1)N\big)+2\big)e^{-y}\leq 2e\big(y\big[\ln((y-1)N\big)+2\big]+2\big)e^{-y}.
$$
Next, we observe that  $\bar s_N\leq \bar \mu_N(1-\bar\mu_N)$, and replacing $\bar s_N$ in \rf{mart10} with this upper bound come to the inequality:
\[
\Prob\left\{|\bar \zeta_N|\geq \sqrt{2y\bar \mu_N(1-\bar\mu_N)\over N }+{y\over 3N }\right\}\leq 2e \big(y\big[\ln((y-1)N\big)+2\big]+2\big) e^{-y}.
\]
In other words, there exist a subset $\overline \Omega^N$ of the space $\Omega^N$ of realizations $\omega^N$ of probability at least $1-2e\big(y\ln\big((y-1)n\big)+4\big)e^{-y}$ and such for all $\omega^N\in \overline \Omega^N$ one has
\be
%\overline\phi(\mu,I;y)\leq \bar \zeta_n\;\;\mbox{and}\;\;\underline\phi(\mu,I;y)\geq \bar \zeta_n
%-\sqrt{2y \bar\mu_N(1-\bar\mu_n)\over N}-{y\over 3N}\leq
|\bar\zeta_N|\leq \sqrt{2y \bar\mu_N(1-\bar\mu_N)\over N}+{y\over 3N}.
\ee{simplea}
%(note that the latter constraints are convex in $\bar \mu_n$).
Observe that $\bar\mu_n$ can be eliminated from the above inequalities: when denoting
$\nu_i=\gamma_i^T\omega_i%=\zeta_i+\gamma_i^T\mu_i
$ with $\bar\nu_N={1\over N}\sum_{i=1}^N \nu_i=\bar\zeta_N+\bar\mu_N$, by simple algebra we deduce from \rf{simplea} that
\[
\underline\psi(\bar \nu_N,I;y)\leq \bar\mu_N\leq \overline\psi(\bar \nu_N,I;y)
\]
where  $\underline\psi(\cdot)$ and $\overline\psi(\cdot)$ are as in \rf{psifun}.
We conclude that for $\omega^N\in \overline \Omega^N$
$$
\bar \nu_N-\overline\psi(\bar \nu_N,I;y) \leq \bar\zeta_N\leq \bar\nu_N-\underline\psi(\bar \nu_N,I;y)
$$
what implies \rf{mart010}. \qed
\par\noindent{\em Proof of Lemma \ref{lem1002}.}
Now, in the premise of Lemma \ref{lem1002}, let us fix $k\in\{1,...,\kappa\}$, and let us denote
$\gamma_i^T=[\etai]_k=\Row_k[\etai]$, the $k$-th row of $\etai$. We set $\nu_i=\gamma_i^T\omega_i=[\etai]_k\omega_i$. Note that conditional distribution of the r.v. $\nu_i$ given $\omega^{i-1}$ is
Bernoulli distribution with parameter $\mu_i=\bE_{|\omega^{i-1}}\{\nu_i\}=[\etai]_k\etait\beta$. Defining, as above, $\zeta_i=\nu_i-\mu_i$, $\bar\zeta_N={1\over N}\sum_{i=1}^N \zeta_i=F_{\omega^N}(\beta)_k$, the $k$-th component of the field $F_{\omega^N}(\beta)$, $\bar\nu_N={1\over N}\sum_{i=1}^N \nu_i=a[\omega^N]_k$, and
$\bar\mu_N={1\over N}\sum_{i=1}^N \mu_i= {1\over N}\sum_{i=1}^N  [\etai]_k\etait\beta=(A[\omega^N]\beta)_k$, the $k$-th component of $A[\omega^N]\beta$,
and utilizing bound \rf{mart010} of Lemma \ref{lem:mart} we conclude that for any $y>1$ $F_{\omega^N}(\beta)_k$, $k=1,...,\kappa$,
satisfy, with probability at least $1-2e \big(y\big[\ln((y-1)N\big)+2\big]+2\big) e^{-y}$, the bound
\[
\bar \nu_N-\overline\psi(\bar \nu_N,N;y)\leq F_{\omega^N}(\beta)_k\leq \bar\nu_N-\underline\psi(\bar \nu_N,N;y)
\]
where $\underline\psi(\cdot)$ and $\overline\psi(\cdot)$ are as in \rf{psifun}.\qed

\end{document}